\documentclass[11pt,leqno]{article}
\usepackage{amsfonts}
\usepackage{amsmath}
\usepackage{amssymb}
\usepackage{theorem}
\usepackage{latexsym}
\usepackage{mathrsfs}

\newtheorem{theorem}{\bf Theorem}[section]
\newtheorem{lemma}{\bf Lemma}[section]

\newtheorem{remark}{\bf Remark}[section]

\def\thebibliography#1{\section*{References\markboth
 {REFERENCES}{REFERENCES}}\list
 {[\arabic{enumi}]}{\settowidth\labelwidth{[#1]}\leftmargin\labelwidth
 \advance\leftmargin\labelsep
 \usecounter{enumi}}
 \def\newblock{\hskip .11em plus .33em minus .07em}
 \sloppy
 \sfcode`\.=1000\relax}

\begin{document}
\vspace*{0ex}
\begin{center}
{\Large\bf
Isobe--Kakinuma model for water waves \\[0.5ex]
as a higher order shallow water approximation 
}
\end{center}

\begin{center}
Tatsuo Iguchi
\end{center}

\begin{abstract}
We justify rigorously an Isobe--Kakinuma model for water waves as a higher order 
shallow water approximation in the case of a flat bottom. 
It is known that the full water wave equations are approximated by the shallow 
water equations with an error of order $O(\delta^2)$, 
where $\delta$ is a small nondimensional parameter defined as the ratio of the 
mean depth to the typical wavelength. 
The Green--Naghdi equations are known as higher order approximate equations 
to the water wave equations with an error of order $O(\delta^4)$. 
In this paper we show that the Isobe--Kakinuma model is a much higher order 
approximation to the water wave equations with an error of order $O(\delta^6)$. 
\end{abstract}

\section{Introduction}
\label{section:intro}
We are concerned with a mathematically rigorous justification of an Isobe--Kakinuma 
model for the full water wave problem as a higher order shallow water approximation 
in the strongly nonlinear regime in the case of a flat bottom. 
The water wave problem is mathematically formulated as a free boundary problem for 
an irrotational flow of an inviscid and incompressible fluid under the gravitational field. 
We consider the water filled in $(n+1)$-dimensional Euclidean space. 
Let $t$ be the time, $x=(x_1,\ldots,x_n)$ the horizontal spatial coordinates, 
and $z$ the vertical spatial coordinate. 
We assume that the water surface and the bottom are represented as $z=\eta(x,t)$ and $z=-h$, 
respectively, where $\eta=\eta(x,t)$ is the surface elevation and $h$ is the mean depth. 
J. C. Luke \cite{Luke1967} showed that the water wave problem has a variational structure 
by giving a Lagrangian in terms of the velocity potential $\Phi=\Phi(x,z,t)$ and the surface 
elevation $\eta$. 
His Lagrangian has the form 
\begin{equation}\label{intro:Luke's Lagrangian}
\mathscr{L}(\Phi,\eta) = \int_{-h}^{\eta(x,t)}\biggl(\partial_t\Phi(x,z,t)
 +\frac12|\nabla_X\Phi(x,z,t)|^2+gz\biggr){\rm d}z
\end{equation}
and the action function is 
\[
\mathscr{J}(\Phi,\eta)
= \int_{t_0}^{t_1}\!\!\!\int_{\Omega}\mathscr{L}(\Phi,\eta){\rm d}x{\rm d}t, 
\]
where $\nabla_X=(\nabla,\partial_z)=(\partial_{x_1},\ldots,\partial_{x_n},\partial_z)$, 
$g$ is the gravitational constant, and $\Omega$ is an appropriate region in $\mathbf{R}^n$. 
J. C. Luke showed that the corresponding Euler--Lagrange equation is exactly the basic 
equations for water waves. 
M. Isobe \cite{Isobe1994, Isobe1994-2} and T. Kakinuma \cite{Kakinuma2000, Kakinuma2001, Kakinuma2003} 
approximated the velocity 
potential in Luke's Lagrangian as 
\[
\Phi(x,z,t) \simeq \sum_{k=0}^K\Psi_k(z)\phi_k(x,t),
\]
where $\{\Psi_k\}$ is an appropriate function system in the vertical coordinate $z$, 
and derived an approximate Lagrangian for $(\eta,\phi_0,\phi_1,\ldots,\phi_K)$. 
The Isobe--Kakinuma model is the corresponding Euler--Lagrange equation for the 
approximated Lagrangian. 
Different choices of the function system $\{\Psi_k\}$ give different Isobe--Kakinuma models. 
In this paper we adopt the approximation 
\begin{equation}\label{intro:approximation}
\Phi(x,z,t) \simeq \phi_0(x,t)+(h+z)^2\phi_1(x,t).
\end{equation}
Plugging this into Luke's Lagrangian \eqref{intro:Luke's Lagrangian} we obtain an 
approximate Lagrangian $\mathscr{L}^{\mbox{\rm\tiny app}}(\phi_0,\phi_1,\eta)$. 
The corresponding Euler--Lagrange equation has the form 
\begin{equation}\label{intro:Isobe-Kakinuma model}
\left\{
 \begin{array}{l}
  \displaystyle
  \partial_t\eta+\nabla\cdot\biggl(H\nabla\phi_0+\frac13H^3\nabla\phi_1\biggr)=0, \\[2ex]
  \displaystyle
  H^2\partial_t\eta+\nabla\cdot\biggl(\frac13H^3\nabla\phi_0+\frac15H^5\nabla\phi_1\biggr) 
   -\frac43H^3\phi_1=0, \\[2ex]
  \displaystyle
  \partial_t\phi_0+H^2\partial_t\phi_1+g\eta 
   +\frac12|\nabla\phi_0|^2+H^2\nabla\phi_0\cdot\nabla\phi_1
   +\frac12H^4|\nabla\phi_1|^2+2H^2(\phi_1)^2=0,
 \end{array}
\right.
\end{equation}
where $H=H(x,t)$ is the depth of the water and is given by $H(x,t)=h+\eta(x,t)$. 
This is the Isobe--Kakinuma model that we are going to consider in this paper. 
For the detailed derivation of this model we refer to Y. Murakami and T. Iguchi \cite{MurakamiIguchi2015}.

In order to compare this Isobe--Kakinuma model with the full water wave problem 
in the shallow water regime, we need to rewrite \eqref{intro:Isobe-Kakinuma model} 
in an appropriate nondimensional form. 
Let $\lambda$ be the typical wave length and introduce a nondimensional parameter 
$\delta$ by the aspect ratio $\delta=h/\lambda$, which measures the shallowness of the water. 
We rescale the independent and the dependent variables by 
\[
  x=\lambda\tilde{x}, \quad z=h\tilde{z}, \quad t=\frac{\lambda}{\sqrt{gh}}\tilde{t}, 
  \quad \phi_0=\lambda\sqrt{gh}\tilde{\phi}_0, \quad 
  \phi_1=\frac{\sqrt{gh}}{\lambda}\tilde{\phi}_1, \quad \eta=h\tilde{\eta}. 
\]
Here we note that these rescaling of dependent variables are related to the 
strongly nonlinear regime of the wave. 
Plugging these into \eqref{intro:Isobe-Kakinuma model} and dropping the tilde sign in the 
notation we obtain 
\begin{equation}\label{intro:IK model}
\left\{
 \begin{array}{l}
  \displaystyle
  \partial_t\eta+\nabla\cdot\biggl(H\nabla\phi_0+\frac13\delta^2H^3\nabla\phi_1\biggr)=0, \\[2ex]
  \displaystyle
  H^2\partial_t\eta+\nabla\cdot\biggl(\frac13H^3\nabla\phi_0+\frac15\delta^2H^5\nabla\phi_1\biggr) 
   -\frac43H^3\phi_1=0, \\[2ex]
  \displaystyle
  \partial_t\phi_0+\delta^2H^2\partial_t\phi_1+\eta \\[1ex]
  \displaystyle\quad
   +\frac12|\nabla\phi_0|^2+\delta^2H^2\nabla\phi_0\cdot\nabla\phi_1
   +\frac12\delta^4H^4|\nabla\phi_1|^2+2\delta^2H^2(\phi_1)^2=0,
 \end{array}
\right.
\end{equation}
where $H(x,t)=1+\eta(x,t)$. 
We consider the initial value problem to this Isobe--Kakinuma model \eqref{intro:IK model} 
under the initial conditions 
\begin{equation}\label{intro:initial conditions}
(\eta,\phi_0,\phi_1)=(\eta_{(0)},\phi_{0(0)},\phi_{1(0)}) \quad\makebox[3em]{at} t=0.
\end{equation}
Unique solvability locally in time of the initial value problem 
\eqref{intro:IK model}--\eqref{intro:initial conditions} and fundamental properties 
of the model, especially, the linear dispersion relation are presented in 
Y. Murakami and T. Iguchi \cite{MurakamiIguchi2015}.

On the other hand, the initial value problem to the full water wave problem 
in Zakharov--Craig--Sulem formulation in the nondimensional form is written as 
\begin{equation}\label{intro:WW}
\left\{
 \begin{array}{l}
  \partial_t\eta-\Lambda(\eta,\delta)\phi=0, \\
  \displaystyle
  \partial_t\phi+\eta+\frac12|\nabla\phi|^2
   -\delta^2\frac{(\Lambda(\eta,\delta)\phi+\nabla\eta\cdot\nabla\phi)^2}{2(1+\delta^2|\nabla\eta|^2)}
   =0,
 \end{array}
\right.
\end{equation}
\begin{equation}\label{intro:IC of WW}
(\eta,\phi)=(\eta_{(0)},\phi_{(0)}),
\end{equation}
where $\phi=\phi(x,t)$ is the trace of the velocity potential $\Phi$ on the water surface 
and $\Lambda(\eta,\delta)$ is the Dirichlet-to-Neumann map for Laplace's equation. 
More precisely, the linear operator $\Lambda(\eta,\delta)$ depending nonlinearly on 
the surface elevation $\eta$ and the parameter $\delta$ is defined by 
\begin{equation}\label{intro:DN}
\Lambda(\eta,\delta)\phi=(\delta^{-2}\partial_z\Phi-\nabla\eta\cdot\nabla\Phi)|_{z=\eta(x,t)},
\end{equation}
where $\Phi$ is a unique solution to the boundary value problem for Laplace's equation 
\begin{equation}\label{intro:BBP}
\left\{
 \begin{array}{lll}
  \delta^2\Delta\Phi+\partial_z^2\Phi=0 & \mbox{in} & -1<z<\eta(x,t), \\
  \Phi=\phi & \mbox{on} & z=\eta(x,t), \\
  \partial_z\Phi=0 & \mbox{on} & z=-1.
 \end{array}
\right.
\end{equation}

It has already been established that the solution to the full water wave problem is 
approximated by the solution to the shallow water equations up to order $O(\delta^2)$, 
that is, we have 
\[
|\eta^{\mbox{\tiny WW}}(x,t)-\eta^{\mbox{\tiny SW}}(x,t)| \lesssim \delta^2
\]
on some time interval independent of $\delta$, where $\eta^{\mbox{\tiny WW}}$ 
and $\eta^{\mbox{\tiny SW}}$ are the solutions to the full water wave and to 
the shallow water equations, respectively. 
For this rigorous justification of the shallow water equations, 
we refer to L. V. Ovsjannikov \cite{Ovsjannikov1974, Ovsjannikov1976} and 
T. Kano and T. Nishida \cite{KanoNishida1979} in the case of analytic initial data and 
Y. A. Li \cite{Li2006}, T. Iguchi \cite{Iguchi2009, Iguchi2011}, and 
B. Alvarez-Samaniego and D. Lannes \cite{Alvarez-SamaniegoLannes2008} 
in the case of initial data in Sobolev spaces. 
Green--Naghdi equations are known as higher order approximate equations to the full 
water wave equations in the shallow water regime, that is, we have 
\[
|\eta^{\mbox{\tiny WW}}(x,t)-\eta^{\mbox{\tiny GN}}(x,t)| \lesssim \delta^4
\]
on some time interval independent of $\delta$, where $\eta^{\mbox{\tiny GN}}$ is a solution to 
the Green--Naghdi equations. 
For this approximation we refer to Y. A. Li \cite{Li2006}, 
B. Alvarez-Samaniego and D. Lannes \cite{Alvarez-SamaniegoLannes2008}, 
and H. Fujiwara and T. Iguchi \cite{FujiwaraIguchi2015}. 
In this paper we will show that the Isobe--Kakinuma model \eqref{intro:IK model} 
is a much higher order approximation to the full water wave equations in the 
shallow water regime, that is, we have 
\begin{equation}\label{intro:error}
|\eta^{\mbox{\tiny WW}}(x,t)-\eta^{\mbox{\tiny IK}}(x,t)| \lesssim \delta^6
\end{equation}
on some time interval independent of $\delta$, where $\eta^{\mbox{\tiny IK}}$ is a solution to 
the Isobe--Kakinuma model \eqref{intro:IK model}. 
Here we remark that Y. Matsuno \cite{Matsuno2015, Matsuno2016} derived extended 
Green--Naghdi equations as higher order shallow water approximations in the strongly 
nonlinear regime. 
His $\delta^{2N}$ model is an approximation of the full water wave equations with an 
error of order $\delta^{2N+2}$. 
Since the linear dispersion relation of his $\delta^4$ model does not have good structures, 
we cannot expect the well-posedness of the initial value problem so that it is hopeless to 
obtain an error estimate of the solutions such as \eqref{intro:error}. 
The linear part of his $\delta^6$ model has a good structure and the solution might 
approximate the solution to the full water wave equations up to order $O(\delta^8)$. 
However, it contains 7th order derivative terms, which are troublesome in a numerical computation. 
Although the Isobe--Kakinuma model \eqref{intro:IK model} is a higher order 
shallow water approximation, it is a system of second order partial differential equations 
and does not contain such higher order derivative terms. 
This is a strong advantage of the Isobe--Kakinuma model.

The contents of this paper are as follows. 
In Section \ref{section:result} we present our main results in this paper, that is, 
uniform estimates of the solution of the initial value problem to the Isobe--Kakinuma model 
\eqref{intro:IK model}--\eqref{intro:initial conditions} on some time interval 
independent of the parameter $\delta$, 
the consistency of the Isobe--Kakinuma model at order $O(\delta^6)$, 
and the rigorous justification of the Isobe--Kakinuma model by establishing 
an error estimate of the solutions such as \eqref{intro:error}. 
In Section \ref{section:strategy} we consider linearized equations of the Isobe--Kakinuma model 
around the rest state in order to explaine a hidden symmetric structure of the model and to 
give an idea to obtain uniform estimates of the solution. 
In Section \ref{section:transformation} we derive a symmetric system of quasilinear equations 
for derivatives of the solution. 
In Section \ref{section:estimate} we give uniform estimates of the solution by using the 
symmetric structure of the model. 
In Section \ref{section:consistency} we show that the Isobe--Kakinuma model is consistent at order 
$(\delta^6)$, that is, the solution to the Isobe--Kakinuma model satisfies the full water wave equations 
with an error of order $O(\delta^6)$. 
In Section \ref{section:justification} we derive an error estimate of the solutions 
by using the stability of the full water wave problem.

\bigskip
\noindent
{\bf Notation}. \ 
We denote by $W^{m,p}(\mathbf{R}^n)$ the $L^p$ Sobolev space of order $m$ on $\mathbf{R}^n$. 
The norms of the Lebesgue space $L^p(\mathbf{R}^n)$ and the Sobolev space $H^m=W^{m,2}(\mathbf{R}^n)$ 
are denoted by $|\cdot|_p$ and $\|\cdot\|_m$, respectively. 
The $L^2$-norm and the $L^2$-inner product are simply denoted by $\|\cdot\|$ and $(\cdot,\cdot)_{L^2}$, 
respectively. 
We put $\partial_t=\frac{\partial}{\partial t}$, $\partial_j=\frac{\partial}{\partial x_j}$, 
and $\partial_z=\frac{\partial}{\partial z}$. 
For a multi-index $\alpha=(\alpha_1,\ldots,\alpha_n)$ we put 
$\partial^{\alpha}=\partial_1^{\alpha_1}\cdots\partial_n^{\alpha_n}$. 
$[P,Q]=PQ-QP$ denotes the commutator.

\section{Main results}
\label{section:result}
\setcounter{equation}{0}
The Isobe--Kakinuma model (\ref{intro:IK model}) is written in the matrix form as 
$$
\left(
 \begin{array}{ccc}
  1 & 0 & 0 \\
  H^2 & 0 & 0 \\
  0 & 1 & \delta^2H^2
 \end{array}
\right)
\partial_t
\left(
 \begin{array}{c}
  \eta \\
  \phi_0 \\
  \phi_1
 \end{array}
\right)
+\{\mbox{spatial derivatives}\} = \mbox{\boldmath$0$}.
$$
Since the coefficient matrix has always the zero eigenvalue, the hypersurface $t=0$ 
in the space-time $\mathbf{R}^n\times\mathbf{R}$ is characteristic for the Isobe--Kakinuma 
model (\ref{intro:IK model}), so that the initial value problem 
(\ref{intro:IK model})--(\ref{intro:initial conditions}) is not solvable in general. 
In fact, if the problem has a solution $(\eta,\phi_0,\phi_1)$, then by eliminating the 
time derivative $\partial_t\eta$ from the first two equations in (\ref{intro:IK model}) 
we see that the solution has to satisfy the relation 
\begin{equation}\label{result:necessary condition 1}
H^2\nabla\cdot\biggl(H\nabla\phi_0+\frac13\delta^2H^3\nabla\phi_1\biggr)
= \nabla\cdot\biggl(\frac13H^3\nabla\phi_0+\frac15\delta^2H^5\nabla\phi_1\biggr) 
   -\frac43H^3\phi_1,
\end{equation}
which is equivalent to 
\begin{equation}\label{result:necessary condition 2}
\frac23\Delta\phi_0+\frac{2}{15}\delta^2H^2\Delta\phi_1+\frac43\phi_1=0.
\end{equation}
Therefore, as a necessary condition the initial data $(\eta_{(0)},\phi_{0(0)},\phi_{1(0)})$ 
have to satisfy this relation for the existence of the solution.

We also need to mention that the initial value problem for the full water wave problem 
\eqref{intro:WW} may be broken unless a so-called generalized Rayleigh--Taylor sign condition 
$-\frac{\partial p}{\partial N} \geq c_0>0$ on the water surface is satisfied, 
where $p$ is the pressure and $N$ is the unit outward normal on the water surface. 
For the Isobe--Kakinuma model \eqref{intro:IK model} the corresponding sign condition 
is written as $a(x,t) \geq c_0>0$, where 
\begin{equation}\label{result:a}
a = 1+2\delta^2H\partial_t\phi_1+2\delta^2H\nabla\phi_0\cdot\nabla\phi_1
 +2\delta^4H^3|\nabla\phi_1|^2+4\delta^2H(\phi_1)^2.
\end{equation}

The following theorem is one of the main results in this paper and asserts the 
existence of the solution to the initial value problem 
\eqref{intro:IK model}--\eqref{intro:initial conditions} with uniform bounds of the 
solution on a time interval independent of the small parameter $\delta$.

\begin{theorem}\label{result:uniform estimate}
Let $M_0, c_0>0$ and $m$ be an integer such that $m>n/2+1$. 
There exist a time $T_1>0$ and constants $C,\delta_1>0$ such that for any 
$\delta\in(0,\delta_1]$ if the initial data $(\eta_{(0)},\phi_{0(0)},\phi_{1(0)})$ 
satisfy the relation \eqref{result:necessary condition 1} and 
\begin{equation}\label{result:condition for ID}
\left\{
 \begin{array}{l}
  \|\eta_{(0)}\|_m+\delta\|\eta_{(0)}\|_{m+1}
   +\|\nabla\phi_{0(0)}\|_m+\delta\|\nabla\phi_{0(0)}\|_{m+1} \\[0.5ex]
  \phantom{\|\eta_{(0)}\|_m}
   +\delta\|\phi_{1(0)}\|_m+\delta^2\|\phi_{1(0)}\|_{m+1}+\delta^3\|\phi_{1(0)}\|_{m+2} \leq M_0, \\[0.5ex]
  1+\eta_{(0)}(x) \geq c_0 \qquad\mbox{for}\quad x\in\mathbf{R}^n,
 \end{array}
\right.
\end{equation}
then the initial value problem \eqref{intro:IK model}--\eqref{intro:initial conditions} 
has a unique solution $(\eta,\phi_0,\phi_1)$ on the time interval $[0,T_1]$. 
Moreover, the solution satisfies the uniform bound: 
\begin{equation}\label{result:uniform bound}
\left\{
 \begin{array}{l}
  \|\eta(t)\|_m+\delta\|\eta(t)\|_{m+1}+\|\nabla\phi_0(t)\|_m+\delta\|\nabla\phi_0(t)\|_{m+1} \\[0.5ex]
  \phantom{\|\eta(t)\|_m}
   +\|\phi_1(t)\|_{m-1}+\delta\|\phi_1(t)\|_m+\delta^2\|\phi_1(t)\|_{m+1}+\delta^3\|\phi_1(t)\|_{m+2} \\[0.5ex]
  \phantom{\|\eta(t)\|_m}
   +\|\partial_t\eta(t)\|_{m-1}+\delta\|\partial_t\eta(t)\|_m
   +\|\partial_t\phi_0(t)\|_m+\delta\|\partial_t\phi_0(t)\|_{m+1} \\[0.5ex]
  \phantom{\|\eta(t)\|_m}
   +\delta\|\partial_t\phi_1(t)\|_{m-1}+\delta^2\|\partial_t\phi_1(t)\|_m
   +\delta^3\|\partial_t\phi_1(t)\|_{m+1} \leq C, \\[0.5ex]
  1+\eta(x,t) \geq c_0/2, \quad a(x,t) \geq 1/2 
   \qquad\mbox{for}\quad x\in\mathbf{R}^n, \; 0\leq t\leq T_1,
 \end{array}
\right.
\end{equation}
where $a$ is defined in terms of the solution by \eqref{result:a}.
\end{theorem}

\begin{remark}
{\rm 
In the above theorem the constant $\delta_1$ is small. 
We can reduce the restriction $0<\delta\leq\delta_1$ to, for example, 
$0<\delta\leq1$, if we impose the sign condition $a(x,0)\geq c_0$ on the initial data. 
However, we are interested in the shallow water approximation, that is, the asymptotic 
behavior of the solution as $\delta\to+0$ so that the condition $0<\delta\leq\delta_1$ 
is not an essential restriction. 
}
\end{remark}

Next, we proceed to show that the water wave equations \eqref{intro:WW} are consistent 
at order $O(\delta^6)$ with the Isobe--Kakinuma model \eqref{intro:IK model}. 
To this end we need to relate the dependent variables $(\eta,\phi_0,\phi_1)$ for \eqref{intro:IK model} 
and $(\eta,\phi)$ for \eqref{intro:WW}. 
In view of the facts that $\phi$ is the trace of the velocity potential $\Phi$ on the water surface 
and that $\phi_0$ and $\phi_1$ appear in the approximation \eqref{intro:approximation}, 
these variables are related by the formula 
\begin{equation}\label{result:relation}
\phi=\phi_0+\delta^2H^2\phi_1
\end{equation}
in the nondimensional variables.

\begin{theorem}\label{result:consistency}
In addition to hypothesis of Theorem {\rm \ref{result:uniform estimate}} we assume that $m>n/2+6$. 
Let $(\eta,\phi_0,\phi_1)$ be the solution obtained in Theorem {\rm \ref{result:uniform estimate}} 
and define $\phi$ by \eqref{result:relation}. 
Then, $(\eta,\phi)$ satisfy the water wave equations with errors of order $O(\delta^6)$, that is, 
\begin{equation}\label{result:approximate WW}
\left\{
 \begin{array}{l}
  \partial_t\eta-\Lambda(\eta,\delta)\phi=\delta^6r_1, \\
  \displaystyle
  \partial_t\phi+\eta+\frac12|\nabla\phi|^2
   -\delta^2\frac{(\Lambda(\eta,\delta)\phi+\nabla\eta\cdot\nabla\phi)^2}{2(1+\delta^2|\nabla\eta|^2)}
   =\delta^6r_2.
 \end{array}
\right.
\end{equation}
Here, $(r_1,r_2)$ satisfy the uniform bound: 
\begin{equation}\label{result:error estimate}
\|r_1(t)\|_{m-7}+\|r_2(t)\|_{m-5} \leq C 
 \qquad\mbox{for}\quad 0\leq t\leq T_1,
\end{equation}
where $C$ is a positive constant independent of $\delta\in(0,\delta_1]$. 
\end{theorem}

The above theorem concerns the approximation of the equations. 
Next, we will be concerned with the approximation of the solution to give a 
rigorous justification of the Isobe--Kakinuma model. 
Here we recall the existence theorem for the initial value problem to the full water wave equations 
\eqref{intro:WW}--\eqref{intro:IC of WW} obtained by T. Iguchi \cite{Iguchi2009}. 
Similar results are obtained by B. Alvarez-Samaniego and D. Lannes 
\cite{Alvarez-SamaniegoLannes2008} and D. Lannes \cite{Lannes2013-2}.

\begin{theorem}\label{result:existence WW}
Let $M_0, c_0>0$ and $m>n/2+1$. 
There exist a time $T_2>0$ and constants $C,\delta_2>0$ such that for any 
$\delta\in(0,\delta_2]$ if the initial data $(\eta_{(0)},\phi_{(0)})$ satisfy 
\[
\left\{
 \begin{array}{l}
  \|\eta_{(0)}\|_{m+3+1/2}+\|\nabla\phi_{(0)}\|_{m+3} \leq M_0, \\[0.5ex]
  1+\eta_{(0)}(x) \geq c_0 \qquad\mbox{for}\quad x\in\mathbf{R}^n,
 \end{array}
\right.
\]
then the initial value problem \eqref{intro:WW}--\eqref{intro:IC of WW} 
has a unique solution $(\eta,\phi)$ on the time interval $[0,T_2]$. 
Moreover, the solution satisfies the uniform bound: 
\[
\left\{
 \begin{array}{l}
  \|\eta(t)\|_{m+3}+\|\nabla\phi(t)\|_{m+2}
   +\|\partial_t\eta(t)\|_{m+2}+\|\partial_t\phi(t)\|_{m+2} \leq C, \\[0.5ex]
  1+\eta(x,t) \geq c_0/2, 
   \qquad\mbox{for}\quad x\in\mathbf{R}^n, \; 0\leq t\leq T_2. 
 \end{array}
\right.
\]
\end{theorem}

\begin{remark}
{\rm 
In the above theorem the constant $\delta_2$ is small. 
As in the case of Theorem \ref{result:uniform estimate} we can reduce the restriction 
$0<\delta\leq\delta_2$ to $0\leq\delta\leq1$, if we impose the sign condition 
$a^{\mbox{\rm\tiny WW}}(x,0)\geq c_0$ on the initial data, where 
$a^{\mbox{\rm\tiny WW}}=1+\delta^2\partial_tZ+\delta^2\mbox{\boldmath $v$}\cdot\nabla Z$ with 
\[
\left\{
 \begin{array}{l}
  Z=(1+\delta^2|\nabla\eta|^2)^{-1}(\Lambda(\eta,\delta)\phi+\nabla\eta\cdot\nabla\phi), \\[0.5ex]
  \mbox{\boldmath $v$}=\nabla\phi-\delta^2Z\nabla\eta.
 \end{array}
\right.
\]

}
\end{remark}

In order that the solution to the Isobe--Kakinuma model 
\eqref{intro:IK model}--\eqref{intro:initial conditions} approximates the solution to the 
full water wave problem \eqref{intro:WW}--\eqref{intro:IC of WW}, 
we need to prepare the initial data $\phi_{0(0)}$ and $\phi_{1(0)}$ for the Isobe--Kakinuma model 
appropriately in terms of the initial data $\eta_{(0)}$ and $\phi_{(0)}$ for the water wave problem. 
In view of the necessary condition \eqref{result:necessary condition 1} or 
\eqref{result:necessary condition 2} and the relation \eqref{result:relation}, 
the initial data have to satisfy 
\begin{equation}\label{result:ID for IK}
\left\{
 \begin{array}{l}
  \displaystyle
  \frac12\Delta\phi_{0(0)}+\frac{1}{10}\delta^2(1+\eta_{(0)})^2\Delta\phi_{1(0)}+\phi_{1(0)}=0, \\[1.5ex]
  \phi_{0(0)}+\delta^2(1+\eta_{(0)})^2\phi_{1(0)}=\phi_{(0)}.
 \end{array}
\right.
\end{equation}
As we will see in Section \ref{section:justification} (see also Lemma \ref{estimate:elliptic 3} and 
Remark \ref{estimate:remark}), 
given the initial data $(\eta_{(0)},\phi_{(0)})$, 
these equations determine uniquely the initial data $(\phi_{0(0)},\phi_{1(0)})$. 
The next theorem gives a rigorous justification of the Isobe--Kakinuma model for the full 
water wave problem as a higher order shallow water approximation.

\begin{theorem}\label{result:justification}
Let $M_0, c_0>0$ and $m$ be an integer such that $m>n/2+1$, and put 
$T_*=\min\{T_1,T_2\}$ and $\delta_*=\min\{\delta_1,\delta_2\}$, where 
these constants $T_1, T_2, \delta_1,\delta_2$ are those in Theorems 
{\rm \ref{result:uniform estimate}} and {\rm \ref{result:existence WW}}. 
Suppose that $0<\delta\leq\delta_*$ and the initial data $(\eta_{(0)},\phi_{(0)})$ satisfy 
\begin{equation}\label{result:condition for ID 2}
\left\{
 \begin{array}{l}
  \|\eta_{(0)}\|_{m+11}+\|\nabla\phi_{(0)}\|_{m+10} \leq M_0, \\[0.5ex]
  1+\eta_{(0)}(x) \geq c_0 \qquad\mbox{for}\quad x\in\mathbf{R}^n. 
 \end{array}
\right.
\end{equation}
Then, \eqref{result:ID for IK} determines uniquely the initial data $(\phi_{0(0)},\phi_{1(0)})$. 
Let $(\eta^{\mbox{\rm\tiny WW}},\phi^{\mbox{\rm\tiny WW}})$ be the solution to the initial value 
problem \eqref{intro:WW}--\eqref{intro:IC of WW} obtained in Theorem {\rm \ref{result:existence WW}} 
and $(\eta^{\mbox{\rm\tiny IK}},\phi_0,\phi_1)$ the solution to the initial value problem 
\eqref{intro:IK model}--\eqref{intro:initial conditions} obtained in Theorem 
{\rm \ref{result:uniform estimate}}, and define $\phi^{\mbox{\rm\tiny IK}}$ by \eqref{result:relation}. 
Then, for any $\delta\in(0,\delta_*]$ we have 
\begin{equation}\label{result:error estimate 2}
\|\eta^{\mbox{\rm\tiny WW}}(t)-\eta^{\mbox{\rm\tiny IK}}(t)\|_{m+2}
 +\|\nabla\phi^{\mbox{\rm\tiny WW}}(t)-\nabla\phi^{\mbox{\rm\tiny IK}}(t)\|_{m+1} \leq C\delta^6
 \qquad\mbox{for}\quad  0\leq t\leq T_*, 
\end{equation}
where $C$ is a positive constant independent of $\delta\in(0,\delta_*]$. 
\end{theorem}

\begin{remark}
{\rm 
The error estimate \eqref{result:error estimate 2} together with the Sobolev imbedding theorem 
implies the pointwise error estimate \eqref{intro:error}. 
}
\end{remark}

We will give the proof of Theorems \ref{result:uniform estimate}, \ref{result:consistency}, and 
\ref{result:justification} in Sections \ref{section:estimate}, \ref{section:consistency}, and 
\ref{section:justification}, respectively.

\section{Strategy to obtain uniform estimates}
\label{section:strategy}
\setcounter{equation}{0}
The unique existence of the solution locally in time to the initial value problem for the 
Isobe--Kakinuma model \eqref{intro:IK model}--\eqref{intro:initial conditions} for each 
fixed $\delta>0$ is established in Y. Murakami and T. Iguchi \cite{MurakamiIguchi2015}. 
However, the energy method used in \cite{MurakamiIguchi2015} does not give uniform estimates 
of the solution with respect to the small parameter $\delta$. 
In order to obtain such estimates we have to make use of a good symmetric structure of 
the Isobe--Kakinuma model \eqref{intro:IK model}. 
In this section we treat linearized equations of the Isobe--Kakinuma model to give an idea 
for obtaining such estimates.

We note that $(\eta,\phi_0,\phi_1)=\mbox{\boldmath$0$}$ is the solution of the Isobe--Kakinuma model 
\eqref{intro:IK model}, which corresponds to the still water with flat water surface. 
The linearized equations of \eqref{intro:IK model} around this trivial solution have the form 
\begin{equation}\label{strategy:linearized equations}
\left\{
 \begin{array}{l}
  \displaystyle
  \partial_t\eta+\Delta\phi_0+\frac13\delta^2\Delta\phi_1=0, \\[1.5ex]
  \displaystyle
  \partial_t\eta+\frac13\Delta\phi_0+\frac15\delta^2\Delta\phi_1-\frac43\phi_1=0, \\[1.5ex]
  \displaystyle
  \partial_t\phi_0+\delta^2\partial_t\phi_1+\eta=0.
 \end{array}
\right.
\end{equation}
Put $\mbox{\boldmath $U$}=(\eta,\phi_0,\phi_1)^{\rm T}$ and 
\[
S=
\left(
 \begin{array}{ccc}
  0         & 1 & \delta^2 \\
  -1        & 0 & 0 \\
  -\delta^2 & 0 & 0
 \end{array}
\right), \quad
A_0(D)=
\left(
 \begin{array}{ccc}
  1 & 0                      & 0 \\
  0 & -\Delta                & -\frac13\delta^2\Delta \\
  0 & -\frac13\delta^2\Delta & -\frac15\delta^4\Delta+\frac43\delta^2
 \end{array}
\right).
\]
Then, the system of equations \eqref{strategy:linearized equations} can be written in 
the matrix form as 
\begin{equation}\label{strategy:matrix form}
S\partial_t\mbox{\boldmath$U$}+A_0(D)\mbox{\boldmath$U$}=\mbox{\boldmath$0$}.
\end{equation}
It is easy to see that $S$ is skew-symmetric and $A_0(D)$ is symmetric in $L^2(\mathbf{R}^n)$. 
Moreover, we can easily show the following lemma, which guarantees the positivity of $A_0(D)$.

\begin{lemma}\label{strategy:positivity}
Let $E(\mbox{\boldmath$U$})=\frac12(A_0(D)\mbox{\boldmath$U$},\mbox{\boldmath$U$})_{L^2}$. 
Then, we have 
\begin{equation}\label{strategy:energy}
E(\mbox{\boldmath$U$}) = \frac12\int_{\mathbf{R}^n}\eta(x)^2\,{\rm d}x
 +\frac12\int_{\mathbf{R}^n}\!\!{\rm d}x\!\int_0^1
  |\nabla_X^{\delta}(\phi_0(x)+\delta^2z^2\phi_1(x))|^2\,{\rm d}z,
\end{equation}
where $\nabla_X^{\delta}=(\nabla,\delta^{-1}\partial_z)$. 
Moreover, $E(\mbox{\boldmath$U$})$ is equivalent to 
\[
\widetilde{E}(\mbox{\boldmath$U$})
=\|\eta\|^2+\|\nabla\phi_0\|^2+\delta^2\|\phi_1\|^2+\delta^4\|\nabla\phi_1\|^2
\]
uniformly with respect to $\delta$. 
\end{lemma}

Let $\mbox{\boldmath$U$}$ be a smooth solution to \eqref{strategy:matrix form}. 
Taking the Euclidean inner product of \eqref{strategy:matrix form} with 
$\partial_t\mbox{\boldmath$U$}$, we obtain 
$A_0(D)\mbox{\boldmath$U$}\cdot\partial_t\mbox{\boldmath$U$}=0$ because $S$ is skew-symmetric. 
Integrating this with respect to $x$ on $\mathbf{R}^n$ we see that 
$\frac{\rm d}{{\rm d}t}E(\mbox{\boldmath$U$}(t))=0$. 
Therefore, $E(\mbox{\boldmath$U$})$ is a conserved quantity for \eqref{strategy:matrix form}. 
In fact, $E(\mbox{\boldmath$U$})$ is the physical energy function: the first term in the 
right-hand side of \eqref{strategy:energy} is the potential energy due to the gravity 
and the second one the kinetic energy. 
However, \eqref{strategy:matrix form} is not standard form of partial differential equations 
because $S$ is a singular matrix. 
In the standard theory of positive systems of partial differential equations, 
the system whose energy function is given by the quadratic form associated to the positive 
operator $A_0(D)$ has the form 
\begin{equation}\label{strategy:positive system}
A_0(D)\partial_t\mbox{\boldmath$U$}+A_1(D)\mbox{\boldmath$U$}=\mbox{\boldmath$0$}
\end{equation}
with a skew-symmetric operator $A_1(D)$ in $L^2(\mathbf{R}^n)$. 
Therefore, we may have a temptation to transform \eqref{strategy:matrix form} into 
\eqref{strategy:positive system}. 
Thus, again let $\mbox{\boldmath$U$}=(\eta,\phi_0,\phi_1)^{\rm T}$ be a smooth solution to 
\eqref{strategy:matrix form}, which is equivalent to \eqref{strategy:linearized equations}, 
and we will derive a system of the form \eqref{strategy:positive system}. 
By eliminating the time derivative $\partial_t\eta$ from the first two equations in 
\eqref{strategy:linearized equations}, we obtain the necessary condition 
\[
\frac12\Delta\phi_0+\frac{1}{10}\delta^2\Delta\phi_1+\phi_1=0
\]
for the existence of the solution. 
We differentiate this with respect to the time $t$. 
The resulting equation together with the third equation in \eqref{strategy:linearized equations} 
can be written in the matrix form 
\[
\left(
 \begin{array}{cc}
  1      & \delta^2 \\
  \frac12\Delta & \frac{1}{10}\delta^2\Delta+1
 \end{array}
\right)\partial_t
\left(
 \begin{array}{c}
  \phi_0 \\
  \phi_1
 \end{array}
\right)+
\left(
 \begin{array}{c}
  \eta \\
  0
 \end{array}
\right)
=\mbox{\boldmath$0$}.
\]
We note that the coefficient matrix operator is invertible, so that this implies 
\begin{equation}\label{strategy:equation for phi}
\left(
 \begin{array}{cc}
  -\Delta                & -\frac13\delta^2\Delta \\
  -\frac13\delta^2\Delta & -\frac15\delta^4\Delta+\frac43\delta^2
 \end{array}
\right)\partial_t
\left(
 \begin{array}{c}
  \phi_0 \\
  \phi_1
 \end{array}
\right)
-c_{IK}(D)^2\Delta
\left(
 \begin{array}{c}
  \eta \\
  \delta^2\eta
 \end{array}
\right)
=\mbox{\boldmath$0$},
\end{equation}
where $c_{IK}(D)^2=(1-\frac25\delta^2\Delta)^{-1}(1-\frac{1}{15}\delta^2\Delta)$. 
We note that the symbol of the operator $c_{IK}(D)$ is the phase speed of the plane 
wave for the linearized Isobe--Kakinuma model. 
This is an evolution equation for $(\phi_0,\phi_1)$. 
We proceed to derive an evolution equation for $\eta$. 
Let $\alpha_j=\alpha_j(D)$ $(j=1,2)$ be Fourier multipliers satisfying $\alpha_1+\alpha_2=1$ 
to be determined later. 
Applying $\alpha_1(D)$ and $\alpha_2(D)$ to the first and the second equations in 
\eqref{strategy:linearized equations}, respectively, we obtain 
\[
\partial_t\eta+\biggl(\Delta\alpha_1+\frac13\Delta\alpha_2\biggr)\phi_0
 +\biggl(\frac13\delta^2\Delta\alpha_1+\biggl(\frac15\delta^2\Delta-\frac43\biggr)\alpha_2\biggr)\phi_1=0.
\]
This and \eqref{strategy:equation for phi} constitute a system of the form 
\eqref{strategy:positive system} with 
\[
A_1(D)=
\left(
 \begin{array}{ccc}
  0 & \Delta\alpha_1+\frac13\Delta\alpha_2
   & \frac13\delta^2\Delta\alpha_1+\bigl(\frac15\delta^2\Delta-\frac43\bigr)\alpha_2 \\
  -c_{IK}(D)^2\Delta & 0 & 0 \\
  -\delta^2c_{IK}(D)^2\Delta & 0 & 0
 \end{array}
\right).
\]
In order that this matrix operator is skew-symmetric in $L^2(\mathbf{R}^n)$, 
the operators $\alpha_1$ and $\alpha_2$ have to satisfy 
\[
\left(
 \begin{array}{cc}
  1 & \frac13 \\
  \frac13\delta^2\Delta & \frac15\delta^2\Delta-\frac43
 \end{array}
\right)
\left(
 \begin{array}{c}
  \alpha_1 \\
  \alpha_2
 \end{array}
\right)
= c_{IK}(D)^2
\left(
 \begin{array}{c}
  1 \\
  \delta^2\Delta
 \end{array}
\right),
\]
which yields 
\[
\left(
 \begin{array}{c}
  \alpha_1 \\
  \alpha_2
 \end{array}
\right)
= c_{IK}(D)^2\biggl(1-\frac{1}{15}\delta^2\Delta\biggr)^{-1}
\left(
 \begin{array}{c}
  1+\frac{1}{10}\delta^2\Delta \\
  -\frac12\delta^2\Delta
 \end{array}
\right).
\]
We note that this choice of $\alpha_1$ and $\alpha_2$ implies the relation 
$\alpha_1+\alpha_2=1$. 
Therefore, we have transformed \eqref{strategy:matrix form} into \eqref{strategy:positive system} 
with $A_1(D)=c_{IK}(D)^2\Delta S$.

However, in this case \eqref{strategy:positive system} is not a system of partial differential 
equations because $c_{IK}(D)$ contains a nonlocal operator $(1-\frac25\delta^2\Delta)^{-1}$. 
Nevertheless, it follows from \eqref{strategy:positive system} that 
\begin{equation}\label{strategy:positive system 2}
\biggl(1-\frac25\delta^2\Delta\biggr)A_0(D)\partial_t\mbox{\boldmath$U$}
 +\biggl(1-\frac{1}{15}\delta^2\Delta\biggr)\Delta S\mbox{\boldmath$U$}=\mbox{\boldmath$0$}.
\end{equation}
This is a positive symmetric system of partial differential equations. 
The corresponding energy function $E_1(\mbox{\boldmath$U$})$ is the quadratic form 
associated with the positive operator $\bigl(1-\frac25\delta^2\Delta\bigr)A_0(D)$, that is, 
\[
E_1(\mbox{\boldmath$U$})
=E(\mbox{\boldmath$U$})+\frac25\delta^2\sum_{j=1}^nE(\partial_{j}\mbox{\boldmath$U$}),
\]
which is equivalent to 
\[
\widetilde{E}_1(\mbox{\boldmath$U$})
= \|\eta\|^2+\delta^2\|\nabla\eta\|^2+\|\nabla\phi_0\|^2+\delta^2\|\Delta\phi_0\|^2
 +\delta^2\|\phi_1\|^2+\delta^4\|\nabla\phi_1\|^2+\delta^6\|\Delta\phi_1\|^2
\]
uniformly with respect to $\delta$. 
Since $(1-\frac{1}{15}\delta^2\Delta)\Delta S$ is skew-symmetric in $L^2(\mathbf{R}^n)$, 
$E_1(\mbox{\boldmath$U$})$ is also a conserved quantity for \eqref{strategy:linearized equations}. 
By using this energy function, 
we can obtain a uniform bound of the solution to \eqref{strategy:linearized equations}.

In view of the above argument, our strategy to obtain uniform estimate of the solution to 
the nonlinear problem is to derive a nonlinear version of the symmetric system 
\eqref{strategy:positive system 2}. 
We will carry out it in the next section.

\section{Transformation of the system}
\label{section:transformation}
\setcounter{equation}{0}
Let $\mbox{\boldmath$U$}=(\eta,\phi_0,\phi_1)^{\rm T}$ be a solution to the Isobe--Kakinuma model 
\eqref{intro:IK model} throughout this section.
We introduce second order differential operators $L_{11}=L_{11}(H)$, $L_{12}=L_{12}(H)$, 
and $L_{22}=L_{22}(H)$ depending on the depth of the water $H=1+\eta$ by 
\begin{equation}\label{transformation:L}
\left\{
 \begin{array}{l}
  L_{11}\psi=-\nabla\cdot(H\nabla\psi), \\[0.5ex]
  \displaystyle
  L_{12}\psi=-\nabla\cdot\biggl(\frac13H^3\nabla\psi\biggr), \\[2ex]
  \displaystyle
  L_{22}\psi=-\delta^2\nabla\cdot\biggl(\frac15H^5\nabla\psi\biggr)+\frac43H^3\psi. 
 \end{array}
\right.
\end{equation}
Then, we see that these operators are symmetric in $L^2(\mathbf{R}^n)$ and that 
the Isobe--Kakinuma model \eqref{intro:IK model} and the relation 
\eqref{result:necessary condition 1} can be written as 
\begin{equation}\label{transformation:IK model}
\left\{
 \begin{array}{l}
  \partial_t\eta-L_{11}\phi_0-\delta^2L_{12}\phi_1=0, \\[0.5ex]
  H^2\partial_t\eta-L_{12}\phi_0-L_{22}\phi_1=0, \\[0.5ex]
  \partial_t\phi_0+\delta^2H^2\partial_t\phi_1+F_1=0
 \end{array}
\right.
\end{equation}
and 
\begin{equation}\label{transformation:necessary condition}
H^2(L_{11}\phi_0+\delta^2L_{12}\phi_1)=L_{12}\phi_0+L_{22}\phi_1,
\end{equation}
respectively, where 
\begin{equation}\label{transformation:F1}
F_1 = \eta+\frac12|\nabla\phi_0|^2+\delta^2H^2\nabla\phi_0\cdot\nabla\phi_1
  +\frac12\delta^4H^4|\nabla\phi_1|^2+2\delta^2H^2(\phi_1)^2. 
\end{equation}
We differentiate \eqref{transformation:necessary condition} (equivalently, 
\eqref{result:necessary condition 2}) with respect to $t$. 
Then, the resulting equation and the third equation in \eqref{transformation:IK model} 
form the system 
\begin{equation}\label{transformation:time detivative}
\left\{
 \begin{array}{l}
  \partial_t\phi_0+\delta^2H^2\partial_t\phi_1=-F_1, \\[0.5ex]
  H^2(L_{11}\partial_t\phi_0+\delta^2L_{12}\partial_t\phi_1)
   =L_{12}\partial_t\phi_0+L_{22}\partial_t\phi_1+F_2,
 \end{array}
\right.
\end{equation}
where 
\begin{equation}\label{transformation:F2}
F_2=\frac{4}{15}\delta^2H^4(\partial_t\eta)\Delta\phi_1.
\end{equation}
Differentiating the equations in \eqref{transformation:time detivative} with respect to 
$t$ once more we obtain 
\begin{equation}\label{transformation:second time detivative}
\left\{
 \begin{array}{l}
  \partial_t^2\phi_0+\delta^2H^2\partial_t^2\phi_1=-F_3, \\[0.5ex]
  H^2(L_{11}\partial_t^2\phi_0+\delta^2L_{12}\partial_t^2\phi_1)
   =L_{12}\partial_t^2\phi_0+L_{22}\partial_t^2\phi_1+F_4,
 \end{array}
\right.
\end{equation}
where 
\begin{equation}\label{transformation:F3F4}
\left\{
 \begin{array}{l}
  F_3=\partial_tF_1+2\delta^2H(\partial_t\eta)(\partial_t\phi_1), \\[1ex]
  \displaystyle
  F_4=\frac{2}{15}\delta^2H^2[\partial_t^2,H^2]\Delta\phi_1.
 \end{array}
\right.
\end{equation}
These systems \eqref{transformation:time detivative} and \eqref{transformation:second time detivative} 
are used to obtain estimates of time derivatives $(\partial_t\phi_0,\partial_t\phi_1)$ and 
$(\partial_t^2\phi_0,\partial_t^2\phi_1)$, respectively.

Let $\alpha=(\alpha_1,\ldots,\alpha_n)$ be a multi-index satisfying $|\alpha|\leq m$. 
We proceed to derive an evolution equation for $\partial^{\alpha}\mbox{\boldmath$U$}$, 
which is a nonlinear version of the symmetric system \eqref{strategy:positive system 2}. 
Applying $\partial^{\alpha}$ to \eqref{transformation:time detivative} we obtain 
\begin{equation}\label{transformation:pre1}
\left(
 \begin{array}{cc}
  1             & \delta^2H^2 \\
  \frac12\Delta & \frac{1}{10}\delta^2H^2\Delta+1
 \end{array}
\right)\partial_t
\left(
 \begin{array}{c}
  \partial^{\alpha}\phi_0 \\
  \partial^{\alpha}\phi_1
 \end{array}
\right)
=
\left(
 \begin{array}{c}
  F_5 \\
  F_6
 \end{array}
\right),
\end{equation}
where 
\begin{equation}\label{transformation:F5F6}
\left\{
 \begin{array}{l}
  F_5=-\delta^2[\partial^{\alpha},H^2]\partial_t\phi_1-\partial^{\alpha}F_1, \\[1ex]
  \displaystyle
  F_6=-\frac{1}{10}\delta^2[\partial^{\alpha},H^2]\Delta\partial_t\phi_1
   -\frac15\delta^2\partial^{\alpha}(H(\partial_t\eta)\Delta\phi_1).
 \end{array}
\right.
\end{equation}
Here, we need to extract principal terms in $F_5$. 
In view of \eqref{transformation:F1} we write $F_1=F_1(\mbox{\boldmath$U$})$ with 
$\mbox{\boldmath$U$}=(\eta,\phi_0,\phi_1)^{\rm T}$. 
We define $a$ by \eqref{result:a} and $\mbox{\boldmath$u$}$ by 
\begin{equation}\label{transformation:u}
\mbox{\boldmath$u$}=\nabla\phi_0+\delta^2H^2\nabla\phi_1,
\end{equation}
which is the horizontal component of the velocity on the water surface. 
Since the Fr\'echet derivative of $F_1(\mbox{\boldmath$U$})$ with respect to 
$\mbox{\boldmath$U$}$ is given by 
\begin{align}\label{transformation:Frechet}
D_{\mbox{\scriptsize\boldmath$U$}}F_1(\mbox{\boldmath$U$})[\zeta,\psi_0,\psi_1]
&= (1+2\delta^2H\nabla\phi_0\cdot\nabla\phi_1
 +2\delta^4H^3|\nabla\phi_1|^2+4\delta^2H(\phi_1)^2)\zeta \\
&\quad
+\mbox{\boldmath$u$}\cdot\nabla\psi_0+\delta^2H^2\mbox{\boldmath$u$}\cdot\nabla\psi_1
 +4\delta^2H^2\phi_1\psi_1, \nonumber
\end{align}
we have 
\begin{equation}\label{transformation:pp of F5}
F_5 = -a\partial^{\alpha}\eta-\mbox{\boldmath$u$}\cdot\nabla\partial^{\alpha}\phi_0
 -\delta^2H^2\mbox{\boldmath$u$}\cdot\nabla\partial^{\alpha}\phi_1+F_7,
\end{equation}
where 
\begin{equation}\label{transformation:F7}
F_7=-\delta^2([\partial^{\alpha},H^2]-2H(\partial^{\alpha}H))\partial_t\phi_1
-(\partial^{\alpha}F_1(\mbox{\boldmath$U$})
 -D_{\mbox{\scriptsize\boldmath$U$}}F_1(\mbox{\boldmath$U$})[\partial^{\alpha}\mbox{\boldmath$U$}])
 -4\delta^2H^2\phi_1\partial^{\alpha}\phi_1.
\end{equation}
Now, we apply the matrix operator 
\[
\left(
 \begin{array}{cc}
  -\Delta(H\,\cdot\,) & -\Delta(\frac13\delta^2H^3\,\cdot\,) \\[0.5ex]
  -\Delta(\frac13\delta^2H^3\,\cdot\,) & -\Delta(\frac15\delta^4H^5\,\cdot\,)+\frac43\delta^2H^3
 \end{array}
\right)
\left(
 \begin{array}{cc}
  \frac{1}{10}\delta^2H^2\Delta+1 & -\delta^2H^2 \\[0.5ex]
  -\frac12\Delta & 1
 \end{array}
\right)
\]
to \eqref{transformation:pre1}. 
In view of the identities 
\begin{align*}
& \left(
 \begin{array}{cc}
  \frac{1}{10}\delta^2H^2\Delta+1 & -\delta^2H^2 \\[0.5ex]
  -\frac12\Delta & 1
 \end{array}
\right)
\left(
 \begin{array}{cc}
  1             & \delta^2H^2 \\[0.5ex]
  \frac12\Delta & \frac{1}{10}\delta^2H^2\Delta+1
 \end{array}
\right) \\
&\quad = 
\left(
 \begin{array}{cc}
  1 & 0 \\
  0 & 1
 \end{array}
\right)
\biggl(1-\frac25\delta^2H^2\Delta\biggr)+
\left(
 \begin{array}{cc}
  0 & \frac{1}{10}\delta^4H^2 \\[0.5ex]
  0 & -\frac12\delta^2
 \end{array}
\right)[\Delta,H^2],
\end{align*}
\begin{align*}
& \left(
 \begin{array}{cc}
  -\Delta(H\,\cdot\,) & -\Delta(\frac13\delta^2H^3\,\cdot\,) \\[0.5ex]
  -\Delta(\frac13\delta^2H^3\,\cdot\,) & -\Delta(\frac15\delta^4H^5\,\cdot\,)+\frac43\delta^2H^3
 \end{array}
\right)
\left(
 \begin{array}{cc}
  0 & \frac{1}{10}\delta^4H^2 \\[0.5ex]
  0 & -\frac12\delta^2
 \end{array}
\right) \\
&\quad = 
\left(
 \begin{array}{cc}
  0 & \Delta(\frac{1}{15}\delta^4H^3\,\cdot\,) \\[0.5ex]
  0 & \Delta(\frac{1}{15}\delta^6H^5\,\cdot\,)-\frac23\delta^4H^3
 \end{array}
\right),
\end{align*}
\begin{align*}
& \left(
 \begin{array}{cc}
  -\Delta(H\,\cdot\,) & -\Delta(\frac13\delta^2H^3\,\cdot\,) \\[0.5ex]
  -\Delta(\frac13\delta^2H^3\,\cdot\,) & -\Delta(\frac15\delta^4H^5\,\cdot\,)+\frac43\delta^2H^3
 \end{array}
\right)
\left(
 \begin{array}{cc}
  \frac{1}{10}\delta^2H^2\Delta+1 & -\delta^2H^2 \\[0.5ex]
  -\frac12\Delta & 1
 \end{array}
\right) \\
&\quad =
\left(
 \begin{array}{cc}
  \Delta(\frac{1}{15}\delta^2H^3\Delta\,\cdot\,)-\Delta(H\,\cdot\,) & 0 \\[0.5ex]
  \Delta(\frac{1}{15}\delta^4H^5\Delta\,\cdot\,)-\Delta(\delta^2H^3\,\cdot\,) & 0
 \end{array}
\right)
+
\left(
 \begin{array}{cc}
  0 & \ \Delta(\frac23\delta^2H^3\,\cdot\,) \\[0.5ex]
  \frac23\delta^2[\Delta,H^3] & \Delta(\frac{2}{15}\delta^4H^5\,\cdot\,)+\frac43\delta^2H^3
 \end{array}
\right),
\end{align*}
\begin{align*}
& \left(
 \begin{array}{cc}
  -\Delta(H\,\cdot\,) & -\Delta(\frac13\delta^2H^3\,\cdot\,) \\[0.5ex]
  -\Delta(\frac13\delta^2H^3\,\cdot\,) & -\Delta(\frac15\delta^4H^5\,\cdot\,)+\frac43\delta^2H^3
 \end{array}
\right)
\biggl(1-\frac25\delta^2H^2\Delta\biggr) \\
&\quad = \mathscr{A}_{22}^{(0)}+
\left(
 \begin{array}{cc}
  -\nabla\cdot((\nabla\eta)\,\cdot\,) & -\nabla\cdot(\delta^2H^2(\nabla\eta)\,\cdot\,) \\[0.5ex]
  -\nabla\cdot(\delta^2H^2(\nabla\eta)\,\cdot\,) & 
   -\nabla\cdot(\delta^4H^4(\nabla\eta)\,\cdot\,)+\frac83\delta^4H^4\nabla\eta\cdot\nabla
 \end{array}
\right),
\end{align*}
\begin{align*}
& \left(
 \begin{array}{c}
  \Delta(\frac{1}{15}\delta^2H^3\Delta\,\cdot\,)-\Delta(H\,\cdot\,) \\[0.5ex]
  \Delta(\frac{1}{15}\delta^4H^5\Delta\,\cdot\,)-\Delta(\delta^2H^3\,\cdot\,) 
 \end{array}
\right)
(\mbox{\boldmath$u$}\cdot\nabla\partial^{\alpha}\phi_0
 +\delta^2H^2\mbox{\boldmath$u$}\cdot\nabla\partial^{\alpha}\phi_1) \\
&\quad = \mathscr{A}_{22}^{(1)}
\left(
 \begin{array}{c}
  \partial^{\alpha}\phi_0 \\
  \partial^{\alpha}\phi_1
 \end{array}
\right)
+
\left(
 \begin{array}{c}
  \Delta(\frac{1}{15}\delta^2H^3[\Delta,\mbox{\boldmath$u$}]\cdot\nabla\partial^{\alpha}\phi_0
   +\frac{1}{15}\delta^4H^3[\Delta,H^2\mbox{\boldmath$u$}]\cdot\nabla\partial^{\alpha}\phi_1) \\[0.5ex]
  \Delta(\frac{1}{15}\delta^4H^5[\Delta,\mbox{\boldmath$u$}]\cdot\nabla\partial^{\alpha}\phi_0
   +\frac{1}{15}\delta^6H^5[\Delta,H^2\mbox{\boldmath$u$}]\cdot\nabla\partial^{\alpha}\phi_1) 
 \end{array}
\right) \\
&\qquad\; -
\left(
 \begin{array}{c}
  \nabla\cdot([\nabla,H(\mbox{\boldmath$u$}\cdot\nabla)]\partial^{\alpha}\phi_0
   +\delta^2[\nabla,H^3(\mbox{\boldmath$u$}\cdot\nabla)]\partial^{\alpha}\phi_1) \\[0.5ex]
  \nabla\cdot(\delta^2[\nabla,H^3(\mbox{\boldmath$u$}\cdot\nabla)]\partial^{\alpha}\phi_0
   +\delta^4[\nabla,H^5(\mbox{\boldmath$u$}\cdot\nabla)]\partial^{\alpha}\phi_1)
 \end{array}
\right),
\end{align*}
where 
\begin{align}\label{transformation:A0}
\mathscr{A}^{(0)}_{22} =&
\left(
 \begin{array}{cc}
  \Delta(\frac25\delta^2H^3\Delta\,\cdot\,)
   & \Delta(\frac{2}{15}\delta^4H^5\Delta\,\cdot\,) \\[0.5ex]
  \Delta(\frac{2}{15}\delta^4H^5\Delta\,\cdot\,)
   & \Delta(\frac{2}{25}\delta^6H^7\Delta\,\cdot\,)
 \end{array}
\right) \\
&- \left(
 \begin{array}{cc}
  \nabla\cdot(H\nabla\,\cdot\,) & \nabla\cdot(\frac13\delta^2H^3\nabla\,\cdot\,) \\[0.5ex]
  \nabla\cdot(\frac13\delta^2H^3\nabla\,\cdot\,) & \nabla\cdot(\frac{11}{15}\delta^4H^5\nabla\,\cdot\,)
 \end{array}
\right)
+
\left(
 \begin{array}{cc}
  0 & 0 \\[0.5ex]
  0 & \frac43\delta^2H^3
 \end{array}
\right), \nonumber
\end{align}
\begin{align}\label{transformation:A1}
\mathscr{A}^{(1)}_{22} =&
\left(
 \begin{array}{cc}
  \Delta(\frac{1}{15}\delta^2H^3(\mbox{\boldmath$u$}\cdot\nabla)\Delta\,\cdot\,)
   & \Delta(\frac{1}{15}\delta^4H^5(\mbox{\boldmath$u$}\cdot\nabla)\Delta\,\cdot\,) \\[0.5ex]
  \Delta(\frac{1}{15}\delta^4H^5(\mbox{\boldmath$u$}\cdot\nabla)\Delta\,\cdot\,)
   & \Delta(\frac{1}{15}\delta^6H^7(\mbox{\boldmath$u$}\cdot\nabla)\Delta\,\cdot\,)
 \end{array}
\right) \\
&- \left(
 \begin{array}{cc}
  \nabla\cdot(H(\mbox{\boldmath$u$}\cdot\nabla)\nabla\,\cdot\,)
   & \nabla\cdot(\delta^2H^3(\mbox{\boldmath$u$}\cdot\nabla)\nabla\,\cdot\,) \\[0.5ex]
  \nabla\cdot(\delta^2H^3(\mbox{\boldmath$u$}\cdot\nabla)\nabla\,\cdot\,)
   & \nabla\cdot(\delta^4H^5(\mbox{\boldmath$u$}\cdot\nabla)\nabla\,\cdot\,)
 \end{array}
\right), \nonumber
\end{align}
we obtain 
\begin{align*}
\mathscr{A}^{(0)}_{22}\partial_t
\left(
 \begin{array}{c}
  \partial^{\alpha}\phi_0 \\
  \partial^{\alpha}\phi_1
 \end{array}
\right)
&+\mathscr{A}^{(1)}_{22}
\left(
 \begin{array}{c}
  \partial^{\alpha}\phi_0 \\
  \partial^{\alpha}\phi_1
 \end{array}
\right)
+
\mathscr{A}^{(1)}_{21}\partial^{\alpha}\eta =
\left(
 \begin{array}{c}
  \delta\Delta G_{1,\alpha}+\nabla\cdot\mbox{\boldmath$G$}_{3,\alpha} \\
  \delta^3\Delta G_{2,\alpha}+\delta^2\nabla\cdot\mbox{\boldmath$G$}_{4,\alpha}+\delta G_{5,\alpha}
 \end{array}
\right),
\end{align*}
where 
\begin{equation}
\mathscr{A}^{(1)}_{21}=
\left(
 \begin{array}{c}
  \Delta(\frac{1}{15}\delta^2H^3\Delta(a\,\cdot\,))-\Delta(Ha\,\cdot\,) \\
  \Delta(\frac{1}{15}\delta^4H^5\Delta(a\,\cdot\,))-\Delta(\delta^2H^3a\,\cdot\,)
 \end{array}
\right)
\end{equation}
and 
\begin{equation}\label{transformation:G1-G4}
\left\{
 \begin{array}{l}
  \displaystyle
  G_{1,\alpha}=\frac{1}{15}\delta^3H^3[\Delta,H^2]\partial_t\partial^{\alpha}\phi_1
   +\frac{1}{15}\delta H^3\Delta F_7+\frac23\delta H^3F_6 \\[2ex]
  \displaystyle
  \phantom{G_{1,\alpha}=}\;
   -\frac{1}{15}\delta H^3[\Delta,\mbox{\boldmath$u$}]\cdot\nabla\partial^{\alpha}\phi_0
   -\frac{1}{15}\delta^3H^3[\Delta,H^2\mbox{\boldmath$u$}]\cdot\nabla\partial^{\alpha}\phi_1, \\[2ex]
  \displaystyle
  G_{2,\alpha}=H^2G_{1,\alpha}-\frac{8}{15}\delta H^5F_6, \\[2ex]
  \displaystyle
  \mbox{\boldmath$G$}_{3,\alpha}
   =(\nabla\eta)(\partial_t\partial^{\alpha}\phi_0+\delta^2H^2\partial_t\partial^{\alpha}\phi_1)
    -\nabla(HF_7) \\[0.5ex]
  \displaystyle
  \phantom{\mbox{\boldmath$G$}_{3,\alpha}=}\;
    +[\nabla,H(\mbox{\boldmath$u$}\cdot\nabla)]\partial^{\alpha}\phi_0
    +\delta^2[\nabla,H^3(\mbox{\boldmath$u$}\cdot\nabla)]\partial^{\alpha}\phi_1, \\[1ex]
  \displaystyle
  \mbox{\boldmath$G$}_{4,\alpha}
   =H^2(\nabla\eta)(\partial_t\partial^{\alpha}\phi_0+\delta^2H^2\partial_t\partial^{\alpha}\phi_1)
    -\nabla(H^3F_7) \\[0.5ex]
  \displaystyle
  \phantom{\mbox{\boldmath$G$}_{4,\alpha}=}\;
    +[\nabla,H^3(\mbox{\boldmath$u$}\cdot\nabla)]\partial^{\alpha}\phi_0
    +\delta^2[\nabla,H^5(\mbox{\boldmath$u$}\cdot\nabla)]\partial^{\alpha}\phi_1, \\[1ex]
  \displaystyle
  G_{5,\alpha}=-\frac83\delta^3H^4\nabla\eta\cdot\nabla\partial_t\partial^{\alpha}\phi_1
   -\frac23\delta^3H^3[\Delta,H^2]\partial_t\partial^{\alpha}\phi_1
   +\frac23\delta[\Delta,H^3]F_5+\frac43\delta H^3F_6.
 \end{array}
\right.
\end{equation}
As we will see later, $\mathscr{A}^{(0)}_{22}$ is positive and $\mathscr{A}^{(1)}_{22}$ 
is skew-symmetric in $L^2(\mathbf{R}^n)$ modulo lower order terms. 
This is the evolution equation for $(\partial^{\alpha}\phi_0,\partial^{\alpha}\phi_1)$.

We proceed to derive an evolution equation for $\partial^{\alpha}\eta$. 
In order to obtain the equation which has a good symmetry we need to note that 
\[
\textstyle
(\mathscr{A}^{(1)}_{21})^*
= \bigl(a\Delta(\frac{1}{15}\delta^2H^3\Delta\,\cdot\,)-aH\Delta, \;
 a\Delta(\frac{1}{15}\delta^4H^5\Delta\,\cdot\,)-\delta^2aH^3\Delta\bigr),
\]
where $P^*$ denotes the adjoint operator of $P$ in $L^2(\mathbf{R}^n)$. 
Applying $\partial^{\alpha}$ to the first and the second equations in \eqref{intro:IK model} 
we obtain 
\begin{equation}\label{transformation:pre2}
\left\{
 \begin{array}{l}
  \displaystyle
  \partial_t\partial^{\alpha}\eta+\mbox{\boldmath$u$}\cdot\nabla\partial^{\alpha}\eta
   +H\Delta\partial^{\alpha}\phi_0+\frac13\delta^2H^3\Delta\partial^{\alpha}\phi_1=F_8, \\[1.5ex]
  \displaystyle
  \partial_t\partial^{\alpha}\eta+\mbox{\boldmath$u$}\cdot\nabla\partial^{\alpha}\eta
   +\frac13H\Delta\partial^{\alpha}\phi_0+\frac15\delta^2H^3\Delta\partial^{\alpha}\phi_1
   -\frac43H\partial^{\alpha}\phi_1=F_9,
 \end{array}
\right.
\end{equation}
where 
\begin{equation}\label{transformation:F8F9}
\left\{
 \begin{array}{l}
  \displaystyle
  F_8=-[\partial^{\alpha},\mbox{\boldmath$u$}]\cdot\nabla\eta
   -[\partial^{\alpha},H]\Delta\phi_0-\frac13\delta^2[\partial^{\alpha},H^3]\Delta\phi_1, \\[1.5ex]
  \displaystyle
  F_9=-[\partial^{\alpha},\mbox{\boldmath$u$}]\cdot\nabla\eta
   -\frac13[\partial^{\alpha},H]\Delta\phi_0-\frac15\delta^2[\partial^{\alpha},H^3]\Delta\phi_1
   +\frac43[\partial^{\alpha},H]\phi_1.
 \end{array}
\right.
\end{equation}
Applying the operators 
$a+\frac{1}{10}\delta^2\Delta(aH^2\,\cdot\,)$ and $-\frac12\delta^2\Delta(aH^2\,\cdot\,)$ to 
the first and the second equations in \eqref{transformation:pre2}, respectively, 
adding the resulting equations, and using the 
equality $\Delta(af)=a\Delta f-(\Delta a)f+2\nabla\cdot((\nabla a)f)$, we obtain 
\[
\mathscr{A}^{(0)}_{11}\partial_t\partial^{\alpha}\eta
 +\mathscr{A}^{(1)}_{11}\partial^{\alpha}\eta
 +\mathscr{A}^{(1)}_{12}
 \left(
  \begin{array}{c}
   \partial^{\alpha}\phi_0 \\
   \partial^{\alpha}\phi_1
  \end{array}
 \right)
 = \delta\nabla\cdot\mbox{\boldmath$G$}_{6,\alpha}+G_{7,\alpha}, 
\]
where 
\begin{equation}\label{transformation:A11}
\left\{
 \begin{array}{l}
  \displaystyle
  \mathscr{A}^{(0)}_{11}=a-\nabla\cdot\biggl(\frac25\delta^2aH^2\nabla\,\cdot\,\biggr), \\[1.5ex]
  \displaystyle
  \mathscr{A}^{(1)}_{11}=a(\mbox{\boldmath$u$}\cdot\nabla)
   -\nabla\cdot\biggl(\frac25\delta^2aH^2(\mbox{\boldmath$u$}\cdot\nabla)\nabla\,\cdot\,\biggr), \\[0.5ex]
  \mathscr{A}^{(1)}_{12}=-(\mathscr{A}^{(1)}_{21})^*
 \end{array}
\right.
\end{equation}
and 
\begin{equation}\label{transformation:G5G6}
\left\{
 \begin{array}{l}
  \displaystyle
  \mbox{\boldmath$G$}_{6,\alpha}=\frac25\delta[\nabla,aH^2]\partial_t\partial^{\alpha}\eta
   +\frac25\delta[\nabla,aH^2(\mbox{\boldmath$u$}\cdot\nabla)]\partial^{\alpha}\eta \\[1.5ex]
  \displaystyle
  \phantom{\mbox{\boldmath$G$}_6=}\;
   +\frac{2}{15}(\nabla a)(\delta H^3\Delta\partial^{\alpha}\phi_0
    +\delta^3H^5\Delta\partial^{\alpha}\phi_1)
   +\frac{1}{10}\delta\nabla(aH^2F_8)-\frac12\delta\nabla(aH^2F_9),\\[2ex]
  \displaystyle
  G_{7,\alpha}=-\frac{1}{15}(\Delta a)(\delta^2H^3\Delta\partial^{\alpha}\phi_0
    +\delta^4H^5\Delta\partial^{\alpha}\phi_1)
   -\frac23\delta^2[\Delta,aH^3]\partial^{\alpha}\phi_1+aF_8.
 \end{array}
\right.
\end{equation}
This is the evolution equation for $\partial^{\alpha}\eta$.

To summarize we derived the evolution equations for $\partial^{\alpha}\mbox{\boldmath$U$}$:
\begin{equation}\label{transformation:quasilinear system}
\mathscr{A}^{(0)}\partial_t\partial^{\alpha}\mbox{\boldmath$U$}
+\mathscr{A}^{(1)}\partial^{\alpha}\mbox{\boldmath$U$} = \mbox{\boldmath$G$}_{\alpha},
\end{equation}
where 
\begin{equation}\label{transformation:A}
\mathscr{A}^{(0)}=
\left(
 \begin{array}{cc}
  \mathscr{A}^{(0)}_{11} & 0 \\
  0 & \mathscr{A}^{(0)}_{22}
 \end{array}
\right), \quad
\mathscr{A}^{(1)}=
\left(
 \begin{array}{cc}
  \mathscr{A}^{(1)}_{11} & \mathscr{A}^{(1)}_{12} \\
  \mathscr{A}^{(1)}_{21} & \mathscr{A}^{(1)}_{22}
 \end{array}
\right),
\end{equation}
and 
\begin{equation}
\mbox{\boldmath$G$}_{\alpha}=
\left(
 \begin{array}{c}
  \delta\nabla\cdot\mbox{\boldmath$G$}_{6,\alpha}+G_{7,\alpha} \\
  \delta\Delta G_{1,\alpha}+\nabla\cdot\mbox{\boldmath$G$}_{3,\alpha} \\
  \delta^3\Delta G_{2,\alpha}+\delta^2\nabla\cdot\mbox{\boldmath$G$}_{4,\alpha}+\delta G_{5,\alpha}
 \end{array}
\right).
\end{equation}
Using these equations we will derive uniform bounds of the solution in the next section.

\section{Uniform estimates}
\label{section:estimate}
\setcounter{equation}{0}
In this section we will prove Theorem \ref{result:uniform estimate}. 
Since the existence theorem has already been established in \cite{MurakamiIguchi2015}, 
it is sufficient to give a priori estimates of the solution. 
In the following of this paper we assume that $0<\delta\leq1$.

In view of the equations \eqref{transformation:time detivative} and 
\eqref{transformation:second time detivative} for the time derivatives and 
\eqref{result:ID for IK} for the initial data, we consider the following elliptic 
partial differential equations for $(\psi_0,\psi_1)$: 
\begin{equation}\label{estimate:elliptic}
\left\{
 \begin{array}{l}
  \psi_0+\delta^2H^2\psi_1=f_1, \\
  H^2(L_{11}\psi_0+\delta^2L_{12}\psi_1)=L_{12}\psi_0+L_{22}\psi_1+f_2+\nabla\cdot\mbox{\boldmath$f$}_3,
 \end{array}
\right.
\end{equation}
where $H=1+\eta$ and the operators $L_{11},L_{12},L_{22}$ are those in 
\eqref{transformation:L}. 
It follows from the first equation in \eqref{estimate:elliptic} that 
$\psi_0=f_1-\delta^2H^2\psi_1$. 
Plugging this into the second equation in \eqref{estimate:elliptic} to eliminate $\psi_0$, 
we obtain 
\begin{equation}\label{estimate:elliptic 2}
L_1\psi_1=-\nabla\cdot\biggl(\frac23H^3\nabla f_1+\mbox{\boldmath$f$}_3\biggr)
 +2H^2\nabla\eta\cdot\nabla f_1-f_2,
\end{equation}
where $L_1=L_1(H)$ is a second order differential operator defined by 
\begin{equation}\label{estimate:L}
L_1\psi_1=\delta^2(H^2L_{11}-L_{12})(H^2\psi_1)+(L_{22}-\delta^2H^2L_{12})\psi_1. 
\end{equation}
We note that the operator $L_1$ is symmetric in $L^2(\mathbf{R}^2)$. 
As was shown in \cite{MurakamiIguchi2015} that the operator $L_1$ is positive in $L^2(\mathbf{R}^n)$. 
More precisely, we have the following lemma.

\begin{lemma}\label{estimate:coercive}
Suppose that $H(x)\geq c_0>0$. 
There exists a positive constant $C=C(c_0)$ depending only on $c_0$ such that we have 
\[
(L_1\psi_1,\psi_1)_{L^2} \geq C^{-1}(\|\psi_1\|^2+\delta^2\|\nabla\psi_1\|^2).
\]
\end{lemma}

\noindent
{\bf Proof}. \ 
We can prove the above estimate in exactly the same way as in \cite{MurakamiIguchi2015}. 
For the sake of completeness, we sketch the proof. 
By direct calculation we have 
\begin{align*}
& (L_{11}\psi_0+\delta^2L_{12}\psi_1,\psi_0)_{L^2}
 +(\delta^2L_{12}\psi_0+\delta^2L_{22}\psi_1,\psi_1)_{L^2} \\
&= \int_{\mathbf{R}^n}\!\!{\rm d}x\!\int_0^{H(x)}
 \bigl(|\nabla\psi_0(x)+\delta^2z^2\nabla\psi_1(x)|^2+(2\delta z\psi_1(x))^2\bigr)\,{\rm d}z \\
&\geq C^{-1} (\|\nabla\psi_0\|^2+\delta^2\|\psi_1\|^2+\delta^4\|\nabla\psi_1\|^2) \\
&\geq C^{-1} (\delta^2\|\psi_1\|^2+\delta^4\|\nabla\psi_1\|^2).
\end{align*}
Therefore, by the definition \eqref{estimate:L} of the operator $L_1$ we see that 
\begin{align*}
(L_1\psi_1,\psi_1)_{L^2}
&= (L_{11}(-\delta H^2\psi_1)+\delta^2L_{12}(\delta^{-1}\psi_1),(-\delta H^2\psi_1))_{L^2} \\
&\quad +(\delta^2L_{12}(-\delta H^2\psi_1)+\delta^2 L_{22}(\delta^{-1}\psi_1),\delta^{-1}\psi_1)_{L^2} \\
&\geq C^{-1}(\|\psi_1\|^2+\delta^2\|\nabla\psi_1\|^2).
\end{align*}
This gives the desired estimate. 
\quad$\Box$

\bigskip
Once we obtain this type of estimate, we can easily show the unique existence of the solution 
to \eqref{estimate:elliptic 2} so that to \eqref{estimate:elliptic} in an appropriate Sobolev space 
by using the standard theory of elliptic partial differential equations. 
Concerning uniform estimates of the solution with respect to $\delta$ we have the following lemma.

\begin{lemma}\label{estimate:elliptic 3}
Let $M_0,c_0>0$ and $m$ be an integer such that $m>n/2+1$. 
There exists a positive constant $C$ such that if $\eta$ and $\delta\in(0,1]$ satisfy 
\[
\|\eta\|_m+\delta\|\eta\|_{m+1} \leq M_0, \quad
c_0\leq H(x)=1+\eta(x) \makebox[3em]{for} x\in\mathbf{R}^n,
\]
then for any $\nabla f_1,f_2,\mbox{\boldmath$f$}_3\in H^l$ with $0\leq l\leq m$, 
\eqref{estimate:elliptic} has a unique solution $(\psi_0,\psi_1)$ satisfying 
\begin{equation}\label{estimate:estimate 1}
\|\nabla\psi_0\|_l^2+\delta^2\|\psi_1\|_l^2+\delta^4\|\nabla\psi_1\|_l^2
\leq C(\|\nabla f_1\|_l^2+\|\mbox{\boldmath$f$}_3\|_l^2+\delta^2\|f_2\|_l^2).
\end{equation}
If in addition $f_1\in H^{l+1}$, then we have 
\begin{equation}\label{estimate:estimate 2}
\|\psi_0\|_l^2 \leq C(\|f_1\|_l^2+\delta^2\|\nabla f_1\|_l^2
 +\delta^2\|\mbox{\boldmath$f$}_3\|_l^2+\delta^4\|f_2\|_l^2).
\end{equation}
\end{lemma}

\begin{remark}\label{estimate:remark}
{\rm 
(1) 
We can reduce the hypothesis $f_2\in H^l$ to $f_2\in H^{l-1}$. 
In this case, the term $\|f_2\|_l^2$ in \eqref{estimate:estimate 1} and 
\eqref{estimate:estimate 2} should be replaced with $\delta^{-2}\|f_2\|_{l-1}^2$. 

(2) 
This lemma guarantees that \eqref{result:ID for IK} determines uniquely the initial data 
$(\eta_{(0)},\phi_{0(0)},\phi_{1(0)})$ for the Isobe--Kakinuma model from the initial data 
$(\eta_{(0)},\phi_{(0)})$ for the full water wave problem. 
}
\end{remark}

\noindent
{\bf Proof}. \ 
Throughout the proof we use the same symbol $C$ to denote positive constants depending only on 
$(M_0,c_0,m)$ and independent of $\delta$. 
It follows from \eqref{estimate:elliptic 2} and Lemma \ref{estimate:coercive} that 
\begin{equation}\label{estimate:pre 1}
\delta^2\|\psi_1\|^2+\delta^4\|\nabla\psi_1\|^2
\leq C(\|\nabla f_1\|^2+\|\mbox{\boldmath$f$}_3\|^2+\delta^2\|f_2\|^2),
\end{equation}
which together with the first equation in \eqref{estimate:elliptic} implies 
\eqref{estimate:estimate 1} in the case $l=0$.

Let $l\geq1$ and $\alpha=(\alpha_1,\ldots,\alpha_n)$ be a multi-index satisfying $|\alpha|\leq l$. 
Applying $\partial^{\alpha}$ to \eqref{estimate:elliptic} we obtain 
\[
\left\{
 \begin{array}{l}
  \displaystyle
  \partial^{\alpha}\psi_0+\delta^2H^2\partial^{\alpha}\psi_1
   =\partial^{\alpha}f_1-\delta^2[\partial^{\alpha},H^2]\psi_1, \\[0.5ex]
  \displaystyle
  H^2(L_{11}\partial^{\alpha}\psi_0+\delta^2L_{12}\partial^{\alpha}\psi_1)
   =L_{12}\partial^{\alpha}\psi_0+L_{22}\partial^{\alpha}\psi_1 \\[0.5ex]
  \displaystyle\qquad\qquad
   +\partial^{\alpha}f_2-\frac{2}{15}\delta^2[\partial^{\alpha},\nabla(H^2)]\cdot\nabla\psi_1
   +\nabla\cdot\biggl(\partial^{\alpha}\mbox{\boldmath$f$}_3
    +\frac{2}{15}\delta^2[\partial^{\alpha},H^2]\nabla\psi_1\biggr). 
 \end{array}
\right.
\]
We apply the estimate obtained just above to $(\partial^{\alpha}\psi_0,\partial^{\alpha}\psi_1)$ 
and obtain 
\begin{equation}\label{estimate:pre 2}
\|\nabla\psi_0\|_l^2+\delta^2\|\psi_1\|_l^2+\delta^4\|\nabla\psi_1\|_l^2
\leq C(\|\nabla f_1\|_l^2+\|\mbox{\boldmath$f$}_3\|_l^2+\delta^2\|f_2\|_l^2+I),
\end{equation}
where 
\[
I=\sum_{|\alpha|\leq l}(\delta^4\|\nabla[\partial^{\alpha},H^2]\psi_1\|^2
 +\delta^4\|[\partial^{\alpha},H^2]\nabla\psi_1\|^2
 +\delta^6\|[\partial^{\alpha},\nabla(H^2)]\cdot\nabla\psi_1\|^2). 
\]
Now, in view of the hypothesis $m>n/2+1$ and $1\leq l\leq m$ we can use the standard commutator estimate 
$\|[\partial^{\alpha},u]v\|\lesssim\|\nabla u\|_{m-1}\|v\|_{l-1}$. 
Since $\nabla[\partial^{\alpha},H^2]\psi_1=[\partial^{\alpha},H^2]\nabla\psi_1
+[\partial^{\alpha},\nabla(H^2)]\psi_1$, we obtain 
\begin{align*}
I &\leq C(\delta^4\|\nabla(H^2)\|_{m-1}^2\|\nabla\psi_1\|_{l-1}^2
 +\delta^4\|\nabla(H^2)\|_m^2\|\psi_1\|_{l-1}^2+\delta^6\|\nabla(H^2)\|_m^2\|\nabla\psi_1\|_{l-1}^2) \\
&\leq C(\delta^4\|\psi_1\|_l^2+\delta^2\|\psi_1\|_{l-1}^2+\delta^4\|\nabla\psi_1\|_{l-1}^2) \\
&\leq \epsilon(\delta^2\|\psi_1\|_l^2+\delta^4\|\nabla\psi_1\|_l^2)
 +C_{\epsilon}\delta^2\|\psi_1\|^2
\end{align*}
for any $\epsilon>0$. 
Plugging this into \eqref{estimate:pre 2}, taking $\epsilon>0$ sufficiently small, and 
using \eqref{estimate:pre 1} we obtain \eqref{estimate:estimate 1}, 
which together with the first equation in \eqref{estimate:elliptic} implies 
\eqref{estimate:estimate 2}. 
\quad$\Box$

\bigskip
Now, we are ready to give a proof of Theorem \ref{result:uniform estimate}. 
Let $\mbox{\boldmath$U$}=(\eta,\phi_0,\phi_1)^{\rm T}$ be a solution of the Isobe--Kakinuma model 
\eqref{intro:IK model}. 
In view of \eqref{transformation:quasilinear system} we define a basic energy function 
$\mathscr{E}(\mbox{\boldmath$U$})=(\mathscr{A}^{(0)}\mbox{\boldmath$U$},\mbox{\boldmath$U$})_{L^2}$. 
It is easy to see that 
\begin{align}\label{estimate:basic energy}
\mathscr{E}(\mbox{\boldmath$U$})
&= (a\eta,\eta)_{L^2}+\frac25\delta^2(aH^2\nabla\eta,\nabla\eta)_{L^2} 
 +(H\nabla\phi_0,\nabla\phi_0)_{L^2} \\
&\quad\;
 +\frac23\delta^2(H^3\nabla\phi_0,\nabla\phi_1)_{L^2}
 +\frac15\delta^4(H^5\nabla\phi_1,\nabla\phi_1)_{L^2}+\frac43\delta^2(H^3\phi_1,\phi_1)_{L^2} \nonumber \\
&\quad\;
 +\frac25\delta^2\biggl\{
  (H^3\Delta\phi_0,\Delta\phi_0)_{L^2}+\frac23\delta^2(H^5\Delta\phi_0,\Delta\phi_1)_{L^2}
  +\frac15\delta^4(H^7\Delta\phi_1,\Delta\phi_1)_{L^2} \nonumber \\
&\makebox[5em]{}
  +\frac43\delta^2(H^5\nabla\phi_1,\nabla\phi_1)_{L^2}\biggr\}, 
 \nonumber
\end{align}
where $H=1+\eta$ and $a$ is the function defined by \eqref{result:a}. 
As is usual, a higher order energy function is defined by 
\[
\mathscr{E}_m(t)=\sum_{|\alpha|\leq m}\mathscr{E}(\partial^{\alpha}\mbox{\boldmath$U$}(t)).
\]
Here we remind that $m$ is assumed to satisfy $m>n/2+1$. 
In view of \eqref{estimate:basic energy} we see that this energy function $\mathscr{E}_m(t)$ 
is equivalent to 
\begin{align*}
E_m(t) 
&= \|\eta(t)\|_m^2+\delta^2\|\eta(t)\|_{m+1}^2+\|\nabla\phi_0(t)\|_m^2
 +\delta^2\|\nabla\phi_0(t)\|_{m+1}^2 \\
&\quad\;
 +\delta^2\|\phi_1(t)\|_m^2+\delta^4\|\phi_1(t)\|_{m+1}^2+\delta^6\|\phi_1(t)\|_{m+2}^2
\end{align*}
uniformly with respect to $\delta\in(0,1]$ under the positivity and the boundedness of 
$H$ and $a$. 
More precisely, we have the following. 
Suppose that the solution $\mbox{\boldmath$U$}=(\eta,\phi_0,\phi_1)^{\rm T}$ satisfies 
\begin{equation}\label{estimate:assumption}
E_m(t)\leq M_1, \quad \frac12c_0\leq H(x,t)\leq 2C_0, \quad \frac12\leq a(x,t) \leq\frac32
\end{equation}
for $x\in\mathbf{R}^n$, $0\leq t\leq T_1$, and $0<\delta\leq \delta_1$, 
where the constants $c_0$ and $C_0$ are determined from the initial datum $\eta_{(0)}$ 
by $c_0\leq H(x,0)\leq C_0$ and 
the constants $M_1$, $T_1$, and $\delta_1$ will be determined later. 
In the following we simply write the constants depending only on 
$(c_0,C_0,m)$ and on $(c_0,C_0,m,M_1)$ by $C_1$ and $C_2$, respectively, 
which may change from line to line. 
Then, it holds that 
\begin{equation}\label{estimate:equivalence}
C_1^{-1}E_m(t) \leq \mathscr{E}_m(t) \leq C_1E_m(t)
\end{equation}
for $0\leq t\leq T_1$ and $0<\delta\leq1$. 
It follows from \eqref{transformation:quasilinear system} that 
\begin{align*}
\frac{\rm d}{{\rm d}t}\mathscr{E}_m(t)
&= \sum_{|\alpha|\leq m}\bigl\{([\partial_t,\mathscr{A}^{(0)}]\partial^{\alpha}\mbox{\boldmath$U$},
  \partial^{\alpha}\mbox{\boldmath$U$})_{L^2}
 +2(\mbox{\boldmath$G$}_{\alpha},\partial^{\alpha}\mbox{\boldmath$U$})_{L^2}
 -2(\mathscr{A}^{(1)}\partial^{\alpha}\mbox{\boldmath$U$},
  \partial^{\alpha}\mbox{\boldmath$U$})_{L^2}\bigr\}.
\end{align*}
Here it is easy to see that 
\[
\left\{
 \begin{array}{l}
  |([\partial_t,\mathscr{A}^{(0)}]\partial^{\alpha}\mbox{\boldmath$U$},
    \partial^{\alpha}\mbox{\boldmath$U$})_{L^2}|
   \leq C_1|(\partial_t\eta,\partial_ta)|_{\infty}\mathscr{E}_m(t), \\[0.5ex]
  |(\mbox{\boldmath$G$}_{\alpha},\partial^{\alpha}\mbox{\boldmath$U$})_{L^2}|
   \leq C_1\mathscr{E}_m(t)+\|(G_{1,\alpha},G_{2,\alpha},\mbox{\boldmath$G$}_{3,\alpha},
    \mbox{\boldmath$G$}_{4,\alpha},G_{5,\alpha},\mbox{\boldmath$G$}_{6,\alpha},
     G_{7,\alpha})\|^2.
 \end{array}
\right.
\]
By the definition \eqref{transformation:A} 
(see also \eqref{transformation:A1} and \eqref{transformation:A11}) and integration by parts, 
we see that 
\begin{align*}
&(\mathscr{A}^{(1)}\partial^{\alpha}\mbox{\boldmath$U$},
  \partial^{\alpha}\mbox{\boldmath$U$})_{L^2} 
= (\mathscr{A}^{(1)}_{11}\partial^{\alpha}\eta,
  \partial^{\alpha}\eta)_{L^2}
 +(\mathscr{A}^{(1)}_{22}
  \left(
   \begin{array}{c}
    \partial^{\alpha}\phi_0 \\
    \partial^{\alpha}\phi_1
   \end{array}
  \right),
  \left(
   \begin{array}{c}
    \partial^{\alpha}\phi_0 \\
    \partial^{\alpha}\phi_1
   \end{array}
  \right))_{L^2} \\
&\quad
= -\frac12((\nabla\cdot(a\mbox{\boldmath$u$}))\partial^{\alpha}\eta,\partial^{\alpha}\eta)_{L^2}
 +\frac15\delta^2((\nabla\cdot(aH^2\mbox{\boldmath$u$}))\nabla\partial^{\alpha}\eta,
  \nabla\partial^{\alpha}\eta)_{L^2} \\
&\qquad\,
 -\frac12((\nabla\cdot(H\mbox{\boldmath$u$}))\nabla\partial^{\alpha}\phi_0,
  \nabla\partial^{\alpha}\phi_0)_{L^2}
 -\delta^2((\nabla\cdot(H^3\mbox{\boldmath$u$}))\nabla\partial^{\alpha}\phi_0,
  \nabla\partial^{\alpha}\phi_1)_{L^2} \\
&\qquad\,
 -\frac12\delta^4((\nabla\cdot(H^5\mbox{\boldmath$u$}))\nabla\partial^{\alpha}\phi_1,
  \nabla\partial^{\alpha}\phi_1)_{L^2}
 -\frac{1}{30}\delta^2((\nabla\cdot(H^3\mbox{\boldmath$u$}))\Delta\partial^{\alpha}\phi_0,
  \Delta\partial^{\alpha}\phi_0)_{L^2} \\
&\qquad\,
 -\frac{1}{15}\delta^4((\nabla\cdot(H^5\mbox{\boldmath$u$}))\Delta\partial^{\alpha}\phi_0,
  \Delta\partial^{\alpha}\phi_1)_{L^2}
 -\frac{1}{30}\delta^6((\nabla\cdot(H^7\mbox{\boldmath$u$}))\Delta\partial^{\alpha}\phi_1,
  \Delta\partial^{\alpha}\phi_1)_{L^2}
\end{align*}
so that 
\[
|(\mathscr{A}^{(1)}\partial^{\alpha}\mbox{\boldmath$U$},
  \partial^{\alpha}\mbox{\boldmath$U$})_{L^2}|
\leq C_1(1+|(\nabla a,\nabla\eta)|_{\infty})\|\mbox{\boldmath$u$}\|_{W^{1,\infty}}
 \mathscr{E}_m(t).
\]
Therefore, we have 
\begin{align}\label{estimate:pre energy estimate}
\frac{\rm d}{{\rm d}t}\mathscr{E}_m(t)
&\leq C_1(1+|(\partial_t\eta,\partial_ta)|_{\infty}
 +(1+|(\nabla\eta,\nabla a)|_{\infty})\|\mbox{\boldmath$u$}\|_{W^{1,\infty}})\mathscr{E}_m(t) \\
&\quad\;
 +\sum_{|\alpha|\leq m}\|(G_{1,\alpha},G_{2,\alpha},\mbox{\boldmath$G$}_{3,\alpha},
    \mbox{\boldmath$G$}_{4,\alpha},G_{5,\alpha},\mbox{\boldmath$G$}_{6,\alpha},
     G_{7,\alpha})\|^2. \nonumber
\end{align}

We proceed to estimate the terms in the right-hand side of the above inequality 
under the condition \eqref{estimate:assumption}. 
In the following, we use the standard calculus inequalities 
\[
\left\{
 \begin{array}{l}
  \displaystyle
  \|uv\|_l \lesssim |u|_{\infty}\|v\|_l+|v|_{\infty}\|u\|_l, \\[0.5ex]
  \|\nabla F(u)\|_l \lesssim C(|u|_{\infty})\|\nabla u\|_l 
 \end{array}
\right.
\]
and the Sobolev imbedding theorem $|u|_{\infty} \lesssim \|u\|_{m-1}$ without any comment. 
It follows from the first equation in \eqref{intro:IK model}, the necessary condition 
\eqref{result:necessary condition 2}, and \eqref{transformation:u} that 
\begin{equation}\label{estimate:pre 3}
\|\partial_t\eta\|_{m-1}^2+\delta^2\|\partial_t\eta\|_m^2+\|\phi_1\|_{m-1}^2
 +\|\mbox{\boldmath$u$}\|_m^2+\delta^2\|\mbox{\boldmath$u$}\|_{m+1}^2
\leq C_2\mathscr{E}_m(t).
\end{equation}
This together with the definitions \eqref{transformation:F1} and \eqref{transformation:F2} 
of $F_1$ and $F_2$ implies 
\[
\|F_1\|_m^2+\delta^2\|F_1\|_{m+1}^2+\|F_2\|_{m-1}^2+\delta^2\|F_2\|_m^2 \leq C_2\mathscr{E}_m(t).
\]
Therefore, applying Lemma \ref{estimate:elliptic 3} with $l=m$ and $l=m-1$ 
to the equation \eqref{transformation:time detivative} 
for $(\partial_t\phi_0,\partial_t\phi_1)$ we obtain 
\begin{equation}\label{estimate:pre 4}
\|\partial_t\phi_0\|_m^2+\delta^2\|\partial_t\phi_0\|_{m+1}^2
+\delta^2\|\partial_t\phi_1\|_{m-1}^2+\delta^4\|\partial_t\phi_1\|_{m}^2
+\delta^6\|\partial_t\phi_1\|_{m+1}^2 \leq C_2\mathscr{E}_m(t).
\end{equation}
Differentiating the first equation in \eqref{intro:IK model} with respect to $t$ we have 
\[
\partial_t^2\eta=-\nabla\cdot\biggl(H\nabla\partial_t\phi_0+\frac13\delta^2H^3\nabla\partial_t\phi_1
 +(\partial_t\eta)\mbox{\boldmath$u$}\biggr)
\]
so that 
\begin{equation}
\delta^2\|\partial_t^2\eta\|_{m-1}^2 \leq C_2\mathscr{E}_m(t).
\end{equation}
In view of \eqref{transformation:Frechet} we have 
\[
\partial_tF_1
= (1+2\delta^2H\mbox{\boldmath$u$}\cdot\nabla\phi_1+4\delta^2H(\phi_1)^2)\partial_t\eta
+\mbox{\boldmath$u$}\cdot\nabla\partial_t\phi_0+\delta^2H^2\mbox{\boldmath$u$}\cdot\nabla\partial_t\phi_1
 +4\delta^2H^2\phi_1\partial_t\phi_1,
\]
which together with the definition \eqref{transformation:F3F4} of $F_3$ and $F_4$ implies 
\[
\delta^2\|F_3\|_m^2+\delta^2\|F_4\|_{m-1} \leq C_2\mathscr{E}_m(t). 
\]
Therefore, applying Lemma \ref{estimate:elliptic 3} with $l=m-1$ to the equation 
\eqref{transformation:second time detivative} 
for $(\partial_t^2\phi_0,\partial_t^2\phi_1)$ we obtain 
\begin{equation}\label{estimate:pre 5}
\delta^4\|\partial_t^2\phi_1\|_{m-1}^2+\delta^6\|\partial_t^2\phi_1\|_m^2
 \leq C_2\mathscr{E}_m(t).
\end{equation}
Then, in view of the definition \eqref{result:a} of $a$ we get 
\begin{equation}\label{estimate:a}
\delta^{-2}\|a-1\|_{m-1}^2+
\|\partial_ta\|_{m-1}^2+\|\nabla a\|_{m-1}^2+\delta^2\|\nabla a\|_m^2 \leq C_2\mathscr{E}_m(t). 
\end{equation}
It follows from \eqref{estimate:assumption}, \eqref{estimate:pre 3}, and \eqref{estimate:a} that 
\begin{equation}\label{estimate:Linfty}
\delta^{-1}|a-1|_{\infty}+
|(\partial_t\eta,\partial_ta)|_{\infty}+\|(\eta,a,\mbox{\boldmath$u$})\|_{W^{1,\infty}}
 +\delta\|(\eta,a,\mbox{\boldmath$u$})\|_{W^{2,\infty}} \leq C_2.
\end{equation}
By the definitions \eqref{transformation:F5F6}, \eqref{transformation:F7}, and 
\eqref{transformation:F8F9} of $F_5,F_6,\ldots,F_9$, we also have 
\begin{equation}\label{estimate:F5-F9}
\|(F_5,F_8,F_9)\|^2+\|F_7\|_1^2
 +\delta^2\bigl(\|F_6\|^2+\|(F_5,F_8,F_9)\|_1^2+\|F_7\|_2^2\bigr) \leq C_2\mathscr{E}_m(t).
\end{equation}
Here, in the estimation of $F_7$ we used the fact that $F_1$ is a polynomial of 
$(\eta,\nabla\phi_0,\delta^2\nabla\phi_1,\delta\phi_1)$ with coefficients independent of $\delta$ 
and the calculus inequality 
\[
\|\partial^{\alpha}F(u)-D_uF(u)[\partial^{\alpha}u]\|_l
 \leq C(\|u\|_{W^{1,\infty}})\|u\|_{|\alpha|+l-1}.
\]
Therefore, by the definitions \eqref{transformation:G1-G4} and \eqref{transformation:G5G6} 
of $G_{j,\alpha}$ we obtain 
\[
\|(G_{1,\alpha},G_{2,\alpha},\mbox{\boldmath$G$}_{3,\alpha},
    \mbox{\boldmath$G$}_{4,\alpha},G_{5,\alpha},\mbox{\boldmath$G$}_{6,\alpha},
     G_{7,\alpha})\|^2
 \leq C_2\mathscr{E}_m(t).
\]
Using this and \eqref{estimate:Linfty} to \eqref{estimate:pre energy estimate} we have 
\[
\frac{\rm d}{{\rm d}t}\mathscr{E}_m(t) \leq C_2\mathscr{E}_m(t),
\]
which together with Gronwall's inequality and \eqref{estimate:equivalence} yields 
$E_m(t)\leq C_1e^{C_2t}M_0^2$, where $M_0$ is the constant in the assumption 
\eqref{result:condition for ID} on the initial data.

To summarize, we have derived the estimates 
\[
E_m(t) \leq C_1e^{C_2t}M_0^2, \quad
 |\eta(x,t)-\eta_{(0)}(x)| \leq C_2t, \quad
 |a(x,t)-1| \leq \delta C_2. 
\]
Taking these into account we define the constants $M_1$, $T_1$, and $\delta_1$ by 
$M_1=2C_1M_0^2$, $T_1=C_2^{-1}\min\{\log2,c_0/2,C_0\}$, and $\delta_1=(2C_2)^{-1}$, 
respectively. 
Then, as in the usual way we can show that the solution actually satisfies \eqref{estimate:assumption}, 
which together with \eqref{estimate:pre 3} and \eqref{estimate:pre 4} yields the uniform estimate 
\eqref{result:uniform bound} of the solution. 
This completes the proof of Theorem \ref{result:uniform estimate}.

\section{Consistency of the Isobe--Kakinuma model}
\label{section:consistency}
\setcounter{equation}{0}
In this section we will prove Theorem \ref{result:consistency}. 
Let $(\eta,\phi_0,\phi_1)$ be a solution of the Isobe--Kakinuma model \eqref{intro:IK model} 
satisfying the uniform bound \eqref{result:uniform bound}, 
and define $\phi$ by the relation \eqref{result:relation}, that is, 
$\phi=\phi_0+\delta^2H^2\phi_1$ with $H=1+\eta$. 
Then, it holds that 
\begin{equation}\label{consistency:estimate phi}
\|\nabla\phi(t)\|_m \leq C \qquad\mbox{for}\quad 0\leq t\leq T_1
\end{equation}
with a constant $C$ independent of $\delta$.

We begin with deriving equations for $(\eta,\phi)$ with errors of order $O(\delta^6)$. 
To this end we need to express $\phi_0$ and $\phi_1$ in terms of $\phi$ and $\eta$. 
Plugging $\phi_0=\phi-\delta^2H^2\phi_1$ into the necessary condition \eqref{result:necessary condition 2} 
we obtain in turn that 
\begin{equation}\label{consistency:phi1}
\left\{
 \begin{array}{l}
  \displaystyle
  \phi_1=-\frac12\Delta\phi+\delta^2R_1, \\[2ex]
  \displaystyle
  \phi_1=-\frac12\Delta\phi
   +\delta^2\biggl\{-\frac14\Delta(H^2\Delta\phi)+\frac{1}{20}H^2\Delta^2\phi\biggr\}
   +\delta^4R_2,
 \end{array}
\right.
\end{equation}
where 
\begin{equation}\label{consistency:R1R2}
\left\{
 \begin{array}{l}
  \displaystyle
  R_1=\frac12\Delta(H^2\phi_1)-\frac{1}{10}H^2\Delta\phi_1, \\[2ex]
  \displaystyle
  R_2=\frac12\Delta(H^2R_1)-\frac{1}{10}H^2\Delta R_1.
 \end{array}
\right.
\end{equation}
Plugging $\phi_0=\phi-\delta^2H^2\phi_1$ into the first equation in \eqref{intro:IK model} 
and using \eqref{consistency:phi1} we obtain 
\begin{equation}\label{consistency:eta}
\left\{
 \begin{array}{l}
  \partial_t\eta=-\nabla\cdot(H\nabla\phi)+\delta^2R_3, \\[0.5ex]
  \displaystyle
  \partial_t\eta=-\nabla\cdot(H\nabla\phi)-\frac13\delta^2\Delta(H^3\Delta\phi)+\delta^4R_4, \\[2ex]
  \displaystyle
  \partial_t\eta=-\nabla\cdot(H\nabla\phi)-\frac13\delta^2\Delta(H^3\Delta\phi) \\[2ex]
  \displaystyle
  \phantom{\partial_t\eta=}
   -\delta^4\biggl\{\frac16\Delta(H^3\Delta(H^2\Delta\phi))-\frac{1}{30}\Delta(H^5\Delta^2\phi)\biggr\}
   +\delta^6R_5,
 \end{array}
\right.
\end{equation}
where 
\begin{equation}\label{consistency:R3-R5}
R_3=\frac23\Delta(H^3\phi_1), \quad R_4=\frac23\Delta(H^3R_1), \quad R_5=\frac23\Delta(H^3R_2). 
\end{equation}
We note that the last equation in \eqref{consistency:eta} can be written in the symmetrical form as 
\begin{equation}\label{consistency:equation for eta}
 \begin{array}{l}
  \displaystyle
  \partial_t\eta+\nabla\cdot(H\nabla\phi)+\frac13\delta^2\Delta(H^3\Delta\phi) \\[1.5ex]
  \displaystyle
  \phantom{\partial_t\eta}
   +\delta^4\biggl\{\frac{1}{15}\Delta(H^3\Delta(H^2\Delta\phi))+\frac{1}{15}\Delta(H^2\Delta(H^3\Delta\phi))
    -\frac15\Delta(|\nabla\eta|^2H^3\Delta\phi)\biggr\}=\delta^6R_5.
 \end{array}
\end{equation}
Plugging $\phi_0=\phi-\delta^2H^2\phi_1$ into the third equation in \eqref{intro:IK model} we obtain 
\[
\partial_t\phi+\eta+\frac12|\nabla\phi|^2
 -2\delta^2H(\phi_1\partial_t\eta+\phi_1\nabla\eta\cdot\nabla\phi)
 +2\delta^2H^2(1+\delta^2|\nabla\eta|^2)(\phi_1)^2=0,
\]
which together with \eqref{consistency:phi1} and \eqref{consistency:eta} yields 
\begin{equation}\label{consistency:equation for phi}
\partial_t\phi+\eta+\frac12|\nabla\phi|^2-\frac12\delta^2H^2(\Delta\phi)^2 
 -\delta^4(\Delta\phi)\biggl\{\frac13H\Delta(H^3\Delta\phi)-\frac12H^2|\nabla\eta|^2\Delta\phi \biggr\}
=\delta^6R_6,
\end{equation}
where 
\begin{align}
R_6 
&= H\biggl\{-(\Delta\phi)R_4-\biggl(\frac12\Delta(H^2\Delta\phi)-\frac{1}{10}H^2\Delta^2\phi\biggr)R_3
 +2(\partial_t\eta+\nabla\eta\cdot\nabla\phi)R_2\biggr\} \\
&\phantom{=}\;
 +H^2\biggl\{(\Delta\phi)R_2+\biggl(\frac12\Delta(H^2\Delta\phi)-\frac{1}{10}H^2\Delta^2\phi\biggr)R_1
 -2\phi_1R_2+|\nabla\eta|^2(\Delta\phi-2\phi_1)R_1\biggr\}. \nonumber
\end{align}
\eqref{consistency:equation for eta} and \eqref{consistency:equation for phi} are the desired 
equations for $(\eta,\phi)$ with errors of order $O(\delta^6)$.

Next, we proceed to expand the full water wave equations with respect to $\delta$ with errors of 
order $O(\delta^6)$. 
To this end we need to expand the Dirichlet-to-Neumann map $\Lambda(\eta,\delta)$ with 
respect to $\delta$. 
It is well-known that $\Lambda(\eta,\delta)$ can be expanded with respect to $\delta^2$ as 
\begin{equation}\label{consistency:expansion}
\Lambda(\eta,\delta)
=\Lambda^{(0)}(\eta)+\delta^2\Lambda^{(1)}(\eta)+\delta^4\Lambda^{(2)}(\eta)+\cdots.
\end{equation}
The explicit forms of these linear operators $\Lambda^{(j)}(\eta)$ are given by the following lemma.

\begin{lemma}\label{DN map}
It holds that 
\[
\left\{
 \begin{array}{l}
  \Lambda^{(0)}(\eta)\psi=-\nabla\cdot(H\nabla\psi), \\[0.5ex]
  \displaystyle
  \Lambda^{(1)}(\eta)\psi=-\frac13\Delta(H^3\Delta\psi), \\[2ex]
  \displaystyle
  \Lambda^{(2)}(\eta)\psi
   =-\frac{1}{15}\Delta(H^3\Delta(H^2\Delta\psi))-\frac{1}{15}\Delta(H^2\Delta(H^3\Delta\psi))
    +\frac15\Delta(|\nabla\eta|^2H^3\Delta\psi).
 \end{array}
\right.
\]
\end{lemma}

\noindent
{\bf Proof.} \ 
These formulae can be derived by the method in T. Iguchi \cite{Iguchi2009} and 
D. Lannes \cite{Lannes2013-2}. 
For the sake of completeness we sketch the proof by following \cite{Lannes2013-2}. 
Let $\Phi$ be the unique solution of the boundary value problem \eqref{intro:BBP}, 
define a diffeomorphism $\Theta$ from the strip $\mathbf{R}^n\times(-1,0)$ to the water region by 
$\Theta(x,z)=(x,z+(z+1)\eta(x))$, and set $\tilde{\Phi}=\Phi\circ\Theta$. 
Then, $\tilde{\Phi}$ satisfies the boundary value problem 
\begin{equation}\label{consistency:BBP}
\left\{
 \begin{array}{lll}
  H^{-1}\partial_z^2\tilde{\Phi}+\delta^2\nabla_X\cdot P\nabla_X\tilde{\Phi}=0
   & \mbox{in} & -1<z<0, \\
  \tilde{\Phi}=\phi & \mbox{on} & z=0, \\
  \partial_z\tilde{\Phi}=0 & \mbox{on} & z=-1,
 \end{array}
\right.
\end{equation}
where $H=1+\eta$, $\nabla_X=(\nabla,\partial_z)$, and 
\[
P(x,z)=
\left(
 \begin{array}{cc}
  H(x)I_n & -(z+1)\nabla\eta(x) \\
  -(z+1)(\nabla\eta(x))^{\rm T} & H^{-1}|\nabla\eta(x)|^2(z+1)^2
 \end{array}
\right)
\]
with the identity matrix $I_n$ of size $n$. 
On the other hand, we define $\overline{\mbox{\boldmath{$V$}}}$ by 
\begin{align*}
\overline{\mbox{\boldmath{$V$}}}(x)
&= \frac{1}{H(x)}\int_{-1}^{\eta(x)}\nabla\Phi(x,z)\,{\rm d}z \\
&= \int_{-1}^0\{\nabla\tilde{\Phi}(x,z)
 -(z+1)H(x)^{-1}(\nabla\eta(x))\partial_z\tilde{\Phi}(x,z)\}\,{\rm d}z,
\end{align*}
which is the vertical average of the horizontal component of the velocity field. 
Then, it holds that 
$\Lambda(\eta,\delta)\phi=-\nabla\cdot(H\overline{\mbox{\boldmath{$V$}}})$. 
Now, expanding $\tilde{\Phi}$ and $\overline{\mbox{\boldmath{$V$}}}$ with respect to $\delta^2$ as 
\[
\left\{
 \begin{array}{l}
  \tilde{\Phi}=\tilde{\Phi}_0+\delta^2\tilde{\Phi}_1+\delta^4\tilde{\Phi}_2+\cdots, \\[0.5ex]
  \overline{\mbox{\boldmath{$V$}}}=\overline{\mbox{\boldmath{$V$}}}_0
   +\delta^2\overline{\mbox{\boldmath{$V$}}}_1+\delta^4\overline{\mbox{\boldmath{$V$}}}_2+\cdots,
 \end{array}
\right.
\]
we have 
\begin{equation}\label{consistency:pre 1}
\left\{
 \begin{array}{l}
  \Lambda^{(j)}(\eta)\phi=-\nabla\cdot(H\overline{\mbox{\boldmath{$V$}}}_j), \\
  \displaystyle
   \overline{\mbox{\boldmath{$V$}}}_j(x)
    = \int_{-1}^0\{\nabla\tilde{\Phi}_j(x,z)
     -(z+1)H(x)^{-1}(\nabla\eta(x))\partial_z\tilde{\Phi}_j(x,z)\}\,{\rm d}z.
 \end{array}
\right.
\end{equation}
Plugging the above expansion of $\tilde{\Phi}$ into \eqref{consistency:BBP} we see that 
$\tilde{\Phi}_0(x,z)=\phi(x)$ and 
\[
\left\{
 \begin{array}{lll}
  H^{-1}\partial_z^2\tilde{\Phi}_j=-\nabla_X\cdot P\nabla_X\tilde{\Phi}_{j-1}
   & \mbox{in} & -1<z<0, \\
  \tilde{\Phi}_j=0 & \mbox{on} & z=0, \\
  \partial_z\tilde{\Phi}_j=0 & \mbox{on} & z=-1
 \end{array}
\right.
\]
for $j=1,2,\ldots$. 
It is not difficult to solve this boundary value problem and we obtain 
\[
\left\{
 \begin{array}{l}
  \displaystyle
  \tilde{\Phi}_1(x,z)=\biggl(-\frac12(z+1)^2+\frac12\biggr)H(x)^2\Delta\phi(x), \\[2ex]
  \displaystyle
  \tilde{\Phi}_2(x,z)=\biggl(\frac18(z+1)^4-\frac14(z+1)^2+\frac18\biggr)
   H(x)^2\Delta(H(x)^2\Delta\phi(x)) \\[2ex]
  \displaystyle
  \phantom{\tilde{\Phi}_2(x,z)=}
   +\biggl(-\frac{1}{12}(z+1)^4+\frac{1}{12}\biggr)H(x)\Delta(H(x)^3\Delta\phi(x)) \\[2ex]
  \displaystyle
  \phantom{\tilde{\Phi}_2(x,z)=}
   +\biggl(\frac14(z+1)^4-\frac14\biggr)H(x)^2|\nabla\eta(x)|^2\Delta\phi(x).
 \end{array}
\right.
\]
Plugging these into \eqref{consistency:pre 1} we obtain the desired formulae. 
\quad$\Box$

\bigskip
By the formulae in this lemma we can rewrite \eqref{consistency:equation for eta} as 
\begin{equation}\label{consistency:equation for eta 2}
\partial_t\eta-\Lambda^{(0)}(\eta)-\delta^2\Lambda^{(1)}(\eta)-\delta^4\Lambda^{(2)}(\eta)
= \delta^6R_5.
\end{equation}
We define remainder terms $R_7,R_8,R_9$ of the expansion \eqref{consistency:expansion} by 
\begin{equation}\label{consistency:R7-R9}
\left\{
 \begin{array}{l}
  \Lambda(\eta,\delta)\phi=\Lambda^{(0)}(\eta)\phi+\delta^2R_7, \\[0.5ex]
  \Lambda(\eta,\delta)\phi
   =\Lambda^{(0)}(\eta)\phi+\delta^2\Lambda^{(1)}(\eta)\phi+\delta^4R_8, \\[0.5ex]
  \Lambda(\eta,\delta)\phi
   =\Lambda^{(0)}(\eta)\phi+\delta^2\Lambda^{(1)}(\eta)\phi+\delta^4\Lambda^{(2)}(\eta)\phi+\delta^6R_9.
 \end{array}
\right.
\end{equation}
In view of the identities 
$(1+\delta^2|\nabla\eta|^2)^{-1}
=1-\delta^2|\nabla\eta|^2+\delta^4|\nabla\eta|^4(1+\delta^2|\nabla\eta|^2)^{-1}$ and 
\begin{align*}
& (\Lambda(\eta,\delta)\phi+\nabla\eta\cdot\nabla\phi)^2 \\
&= (\Lambda^{(0)}(\eta)\phi+\nabla\eta\cdot\nabla\phi)^2
 +\delta^2(\Lambda^{(0)}(\eta)\phi
  +\Lambda(\eta,\delta)\phi+2\nabla\eta\cdot\nabla\phi)R_7 \\
&= (\Lambda^{(0)}(\eta)\phi+\nabla\eta\cdot\nabla\phi)^2
 +2\delta^2(\Lambda^{(0)}(\eta)\phi+\nabla\eta\cdot\nabla\phi)(\Lambda^{(1)}(\eta)\phi) \\
&\quad\;
 +\delta^4\{(\Lambda^{(0)}(\eta)\phi+\Lambda(\eta,\delta)\phi+2\nabla\eta\cdot\nabla\phi)R_8
  +(\Lambda^{(1)}(\eta)\phi)R_7\},
\end{align*}
and Lemma \ref{DN map}, we have 
\begin{align*}
&\frac12\delta^2(1+\delta^2|\nabla\eta|^2)^{-1}(\Lambda(\eta,\delta)\phi+\nabla\eta\cdot\nabla\phi)^2 \\
&= \frac12\delta^2(\Lambda^{(0)}(\eta)\phi+\nabla\eta\cdot\nabla\phi)^2 \\
&\quad
 +\delta^4\biggl\{(\Lambda^{(0)}(\eta)\phi+\nabla\eta\cdot\nabla\phi)(\Lambda^{(1)}(\eta)\phi)
 -\frac12|\nabla\eta|^2(\Lambda^{(0)}(\eta)\phi+\nabla\eta\cdot\nabla\phi)^2\biggr\}
 +\delta^6R_{10} \\
&= \frac12\delta^2H^2(\Delta\phi)^2 
  +\delta^4(\Delta\phi)\biggl\{\frac13H\Delta(H^3\Delta\phi)-\frac12H^2|\nabla\eta|^2\Delta\phi\biggr\}
  +\delta^6R_{10},
\end{align*}
where 
\begin{align}
R_{10}
&= \frac12|\nabla\eta|^4(1+\delta^2|\nabla\eta|^2)^{-1}
 (\Lambda(\eta,\delta)\phi+\nabla\eta\cdot\nabla\phi)^2 \\
&\quad
 -\frac12|\nabla\eta|^2(\Lambda^{(0)}(\eta)\phi
  +\Lambda(\eta,\delta)\phi+2\nabla\eta\cdot\nabla\phi)R_7 \nonumber \\
&\quad
 +\frac12\{(\Lambda^{(0)}(\eta)\phi
  +\Lambda(\eta,\delta)\phi+2\nabla\eta\cdot\nabla\phi)R_8
 +(\Lambda^{(1)}(\eta)\phi)R_7\}. \nonumber
\end{align}
This together with \eqref{consistency:equation for phi}, \eqref{consistency:equation for eta 2}, 
and \eqref{consistency:R7-R9} yields 
\[
\left\{
 \begin{array}{l}
  \partial_t\eta-\Lambda(\eta,\delta)\phi=\delta^6r_1, \\
  \displaystyle
  \partial_t\phi+\eta+\frac12|\nabla\phi|^2
   -\delta^2\frac{(\Lambda(\eta,\delta)\phi+\nabla\eta\cdot\nabla\phi)^2}{2(1+\delta^2|\nabla\eta|^2)}
   =\delta^6r_2,
 \end{array}
\right.
\]
where 
\begin{equation}\label{consistency:r1r2}
r_1=R_5-R_9, \quad r_2=R_6-R_{10}. 
\end{equation}
This is the equation \eqref{result:approximate WW} in Theorem \ref{result:consistency}.

It remains to show the uniform bound \eqref{result:error estimate}. 
To this end we need estimations related to the Dirichlet-to-Neumann map $\Lambda(\eta,\delta)$. 
The following lemma is given in T. Iguchi \cite{Iguchi2009}.

\begin{lemma}\label{consistency:DN}
Let $M_0,c_0>0$ and $l>n/2+1$. 
There exists a positive constant $C$ such that if $\eta\in H^{l+1}$ and $\delta\in(0,1]$ satisfy 
$\|\eta\|_{l+1}\leq M_0$ and $H(x)\geq c_0$ for $x\in\mathbf{R}^n$, then we have 
\[
\|\Lambda(\eta,\delta)\phi\|_l \leq C\|\nabla\phi\|_{l+1}.
\]
\end{lemma}

In order to give systematically error estimates of the expansion \eqref{consistency:expansion} 
it would be better to follow the strategy given by D. Lannes \cite{Lannes2013-2}. 
The following lemma comes easily from the result in \cite{Lannes2013-2}.

\begin{lemma}\label{consistency:error of DN}
Let $M_0,c_0>0$ and suppose that $K$ and $l$ are nonnegative integers such that $l+2K+1>n/2$. 
There exists a positive constant $C$ such that if $\eta\in H^{l+2K+3}$ and $\delta\in(0,1]$ satisfy 
$\|\eta\|_{l+2K+3}\leq M_0$ and $H(x)\geq c_0$ for $x\in\mathbf{R}^n$, then we have 
\[
\|\Lambda(\eta,\delta)\phi-\sum_{k=0}^K\Lambda^{(k)}(\eta)\phi\|_l
 \leq C\delta^{2K+2}\|\nabla\phi\|_{l+2K+3}.
\]
\end{lemma}

In what follows we denote constants depending on $p_1,p_2,\ldots$ by the same symbol 
$C(p_1,p_2,\ldots)$ which may change from line to line. 
Moreover, we may assume that $C(p_1,p_2,\ldots)$ a nondecreasing function of each variable $p_j$. 
Let $l$ be a nonnegative integer and $t_0>n/2$. 
From \eqref{consistency:R1R2} and \eqref{consistency:R3-R5} it follows in turn that 
\[
\left\{
 \begin{array}{l}
  \|(R_1,R_3)\|_l \leq C(\|(\eta,\phi_1)\|_{t_0},\|(\eta,\phi_1)\|_{l+2}), \\[0.5ex]
  \|(R_2,R_4)\|_l \leq C(\|(\eta,R_1)\|_{t_0},\|(\eta,R_1)\|_{l+2}) 
   \leq C(\|(\eta,\phi_1)\|_{t_0+2},\|(\eta,\phi_1)\|_{l+4}), \\[0.5ex]
  \|R_5\|_l \leq C(\|(\eta,R_2)\|_{t_0},\|(\eta,R_2)\|_{l+2}) 
   \leq C(\|(\eta,\phi_1)\|_{t_0+4},\|(\eta,\phi_1)\|_{l+6}),
 \end{array}
\right.
\]
which together with \eqref{consistency:r1r2}, \eqref{consistency:R7-R9}, 
and Lemma \ref{consistency:error of DN} yields 
\[
\|r_1\|_l \leq C(\|(\eta,\phi_1)\|_{t_0+4},\|\phi_1\|_{l+6},\|(\eta,\nabla\phi)\|_{l+7}).
\]
In view of \eqref{result:uniform bound} and \eqref{consistency:estimate phi} choosing 
$t_0=m-5$ and $l=m-7$ we obtain the estimate for $r_1$ in \eqref{result:error estimate}. 
Similarly, we have 
\[
\|R_6\|_l \leq C(\|(\eta,\phi_1)\|_{t_0+4},\|\nabla\phi\|_{t_0+3},\|\partial_t\eta\|_{t_0},
 \|(\eta,\phi_1)\|_{l+4},\|\nabla\phi\|_{l+3},\|\partial_t\eta\|_l).
\]
It follows from \eqref{consistency:R7-R9} and 
Lemmas \ref{consistency:DN} and \ref{consistency:error of DN} that 
\[
\|\Lambda(\eta,\delta)\phi\|_{l+4}+\|R_7\|_{l+2}+\|R_8\|_l
\leq C(\|\eta\|_{l+5})\|\nabla\phi\|_{l+5},
\]
so that 
\[
\|r_2\|_l \leq C(\|(\eta,\nabla\phi)\|_{l+5},\|\phi_1\|_{l+4},\|\partial_t\eta\|_l).
\]
if $l>n/2$. 
By choosing $l=m-5$ we obtain the estimate for $r_2$ in \eqref{result:error estimate}. 
The proof of Theorem \ref{result:consistency} is complete.

\section{Rigorous justification of the Isobe--Kakinuma model}
\label{section:justification}
\setcounter{equation}{0}
In this section we will prove Theorem \ref{result:justification}. 
To this end we take advantage of the stability of the full water wave equations \eqref{intro:WW}, 
which is given by the following theorem. 
Although the statement is not explicitly given in T. Iguchi \cite{Iguchi2009}, 
we can prove it in exactly the same way as the proof of Theorem \ref{result:existence WW}, so that 
we omit the proof. See also D. Lannes \cite{Lannes2013-2}.

\begin{theorem}\label{result:stability}
In addition to hypothesis of Theorem {\rm \ref{result:existence WW}} we assume that 
$0<\delta\leq\delta_1$ and that 
$(\eta^{\mbox{\rm\tiny app}},\phi^{\mbox{\rm\tiny app}})$ satisfy the equations 
\[
\left\{
 \begin{array}{l}
  \partial_t\eta^{\mbox{\rm\tiny app}}
   -\Lambda(\eta^{\mbox{\rm\tiny app}},\delta)\phi^{\mbox{\rm\tiny app}}=f_1^{\mbox{\rm\tiny err}}, \\
  \displaystyle
  \partial_t\phi^{\mbox{\rm\tiny app}}+\eta^{\mbox{\rm\tiny app}}
   +\frac12|\nabla\phi^{\mbox{\rm\tiny app}}|^2
   -\delta^2\frac{(\Lambda(\eta^{\mbox{\rm\tiny app}},\delta)\phi^{\mbox{\rm\tiny app}}
    +\nabla\eta^{\mbox{\rm\tiny app}}\cdot\nabla\phi^{\mbox{\rm\tiny app}})^2}{2(1
     +\delta^2|\nabla\eta^{\mbox{\rm\tiny app}}|^2)}
   =f_2^{\mbox{\rm\tiny err}},
 \end{array}
\right.
\]
on a time interval $[0,T_1]$, the initial condition \eqref{intro:IC of WW}, and the uniform bound: 
\[
\left\{
 \begin{array}{l}
  \|\eta^{\mbox{\rm\tiny app}}(t)\|_{m+3+1/2}+\|\nabla\phi^{\mbox{\rm\tiny app}}(t)\|_{m+3}
   \leq M_1, \\[0.5ex]
  1+\eta^{\mbox{\rm\tiny app}}(x,t) \geq c_0/2
   \qquad\mbox{for}\quad x\in\mathbf{R}^n, \; 0\leq t\leq T_1.
 \end{array}
\right.
\]
Let $(\eta^{\mbox{\rm\tiny WW}},\phi^{\mbox{\rm\tiny WW}})$ be the solution obtained in 
Theorem {\rm \ref{result:existence WW}} and put 
$T_*=\min\{T_1,T_2\}$ and $\delta_*=\min\{\delta_1,\delta_2\}$, where $T_2$ and $\delta_2$ are 
the constants in Theorem {\rm \ref{result:existence WW}}. 
Then, we have 
\begin{align*}
\sup_{0\leq t\leq T_*}
\bigl(\|\eta^{\mbox{\rm\tiny WW}}(t)-\eta^{\mbox{\rm\tiny app}}(t)\|_{m+2}
&  +\|\nabla\phi^{\mbox{\rm\tiny WW}}(t)-\nabla\phi^{\mbox{\rm\tiny app}}(t)\|_{m+1}\bigr) \\
& \leq C_2\sup_{0\leq t\leq T_*}(
 \|f_1^{\mbox{\rm\tiny err}}(t)\|_{m+2}+\|\Lambda_0(\delta)^{1/2}f_2^{\mbox{\rm\tiny err}}(t)\|_{m+2}),
\end{align*}
where $\Lambda_0(\delta)=\Lambda(0,\delta)$ and 
$C_2$ is a positive constant independent of $\delta\in(0,\delta_*]$. 

\end{theorem}

Suppose that the hypotheses in Theorem \ref{result:justification} are satisfied 
for the initial data $(\eta_{(0)},\phi_{(0)})$. 
By Lemma \ref{estimate:elliptic 3} with $m$ replaced by $m+10$ we see that 
\eqref{result:ID for IK} determines uniquely the initial data $(\phi_{0(0)},\phi_{1(0)})$ satisfying 
\[
\|\nabla\phi_{0(0)}\|_{m+10}+\delta\|\phi_{1(0)}\|_{m+10}+\delta^2\|\phi_{1(0)}\|_{m+11} \leq C_0
\]
with a constant $C_0$ independent of $\delta$. 
For these initial data $(\eta_{(0)},\phi_{0(0)},\phi_{1(0)})$ the conditions in Theorems 
\ref{result:uniform estimate} and \ref{result:consistency} with $m$ replaced by $m+9$ are satisfied. 
Therefore, the initial value problem \eqref{intro:IK model}--\eqref{intro:initial conditions} for the 
Isobe--Kakinuma model has a unique solution $(\eta^{\mbox{\rm\tiny IK}},\phi_0,\phi_1)$ 
on the time interval $[0,T_1]$ independent of $\delta\in(0,\delta_1]$. 
Moreover, the solution satisfies the uniform bound \eqref{result:uniform bound} with $m$ replaced 
by $m+9$. 
Put $\phi^{\mbox{\rm\tiny IK}}=\phi_0+\delta^2(1+\eta^{\mbox{\rm\tiny IK}})^2\phi_1$. 
Then, by Theorem \ref{result:consistency} we see that 
$(\eta^{\mbox{\rm\tiny IK}},\phi^{\mbox{\rm\tiny IK}})$ satisfies \eqref{result:approximate WW} with 
$(r_1,r_2)$ satisfying 
\[
\|r_1(t)\|_{m+2}+\|r_2(t)\|_{m+4} \leq C \qquad\mbox{for}\quad 0\leq t\leq T_1,
\]
where $C$ is a constant independent of $\delta\in(0,\delta_1]$. 
Moreover, we have 
\[
\|\eta^{\mbox{\rm\tiny IK}}(t)\|_{m+9}+\|\nabla\phi^{\mbox{\rm\tiny IK}}(t)\|_{m+9}
 \leq C \qquad\mbox{for}\quad 0\leq t\leq T_1.
\]
Therefore, we can apply Theorem \ref{result:stability} and obtain 
\begin{align*}
\sup_{0\leq t\leq T_*}
\bigl(\|\eta^{\mbox{\rm\tiny WW}}(t)
& -\eta^{\mbox{\rm\tiny IK}}(t)\|_{m+2}
 +\|\nabla\phi^{\mbox{\rm\tiny WW}}(t)-\nabla\phi^{\mbox{\rm\tiny IK}}(t)\|_{m+1}\bigr) \\
&\leq C\delta^6\sup_{0\leq t\leq T_*}(
 \|r_1(t)\|_{m+2}+\|\Lambda_0(\delta)^{1/2}r_2(t)\|_{m+2}) \\
&\leq C\delta^6\sup_{0\leq t\leq T_*}(
 \|r_1(t)\|_{m+2}+\|r_2(t)\|_{m+3}) \\
&\leq C\delta^6,
\end{align*}
where we used the estimate $\|\Lambda_0(\delta)^{1/2}\psi\|_s \leq \|\nabla\psi\|_s$. 
This completes the proof of Theorem \ref{result:justification}.


\bigskip
Tatsuo Iguchi \par
{\sc Department of Mathematics} \par
{\sc Faculty of Science and Technology, Keio University} \par
{\sc 3-14-1 Hiyoshi, Kohoku-ku, Yokohama, 223-8522, Japan} \par
E-mail: iguchi@math.keio.ac.jp

\end{document}